\documentclass[12pt]{article}

\usepackage{amsmath,amssymb,amsthm,amscd,mathrsfs}
\usepackage{latexsym,amsmath,amssymb,graphicx}
\usepackage[fleqn]{mathtools}
\usepackage[all]{xy} 
\usepackage{multirow}
\usepackage{color}
\usepackage{epsf}
\usepackage{youngtab}
\xyoption{v2} \xyoption{2cell} \xyoption{dvips}
\CompileMatrices \numberwithin{equation}{section}

\numberwithin{equation}{section}
\makeatletter
\newcommand{\subalign}[1]{%
  \vcenter{%
    \Let@ \restore@math@cr \default@tag
    \baselineskip\fontdimen10 \scriptfont\tw@
    \advance\baselineskip\fontdimen12 \scriptfont\tw@
    \lineskip\thr@@\fontdimen8 \scriptfont\thr@@
    \lineskiplimit\lineskip
    \ialign{\hfil$\m@th\scriptstyle##$&$\m@th\scriptstyle{}##$\crcr
      #1\crcr
    }%
  }
}
\makeatother

\newcommand{\be}{\begin{equation}}
\newcommand{\ee}{\end{equation}}
\newcommand{\IP}{\mathbb{P}}

\newcommand{\oC}{{\overline C}}
\newcommand{\oY}{{\overline Y}}
\newcommand{\oS}{{\overline S}}

\newcommand\IZ{\mathbb {Z}}

\newcommand{\umu}{{\underline \mu}}

\newcommand\IQ{\mathbb {Q}}
\newcommand{\IC}{\mathbb{C}}
\newcommand{\IE}{\mathbb{E}}

\newcommand{\sft}{{\sf t}}

\newcommand{\ba}{\begin{array}}
\newcommand{\ea}{\end{array}}

\newcommand{\wH}{{\widetilde H}}

\newcommand{\om}{\overline{M}}

\newcommand{\unu}{{\underline{\nu}}}
\newcommand{\uell}{{\underline{\ell}}}

\newcommand{\CS}{{\mathcal S}}

\newcommand{\CK}{{\mathcal K}}

\newcommand{\bal}{\begin{aligned}}
\newcommand{\eal}{\end{aligned}}

\newcommand{\sfh}{{\sf u}}

\newcommand{\ualpha}{{\underline \alpha}}

\newcommand{\ch}{{\mathrm{ch}}}

\newcommand{\CO}{{\mathcal O}}

\newcommand{\CH}{{\mathcal H}}

\newcommand{\sfA}{{\sf A}}

\newcommand{\CC}{{\mathcal C}}

\newcommand{\um}{{\underline m}}
\newcommand{\un}{{\underline n}}

\newcommand{\us}{{\underline s}}

\newcommand{\ut}{{\underline t}}
\newcommand{\calP}{{\mathcal P}}

\newcommand{\IA}{{\mathbb A}}
\newcommand{\dl}{{\Delta}}
\newcommand{\uxi}{{\underline \xi}}

\newcommand{\uq}{{\underline q}}

\newcommand{\sfQ}{{\sf Q}}
\newcommand{\sfT}{{\sf T}}
\newcommand{\CY}{{\mathcal Y}}
\CompileMatrices \LaTeXdiagrams \UseAllTwocells
\newdimen\tableauside\tableauside=1.0ex
\newdimen\tableaurule\tableaurule=0.4pt
\newdimen\tableaustep
\def\phantomhrule#1{\hbox{\vbox to0pt{\hrule height\tableaurule width#1\vss}}}
\def\phantomvrule#1{\vbox{\hbox to0pt{\vrule width\tableaurule height#1\hss}}}
\def\sqr{\vbox{%
  \phantomhrule\tableaustep
  \hbox{\phantomvrule\tableaustep\kern\tableaustep\phantomvrule\tableaustep}%
  \hbox{\vbox{\phantomhrule\tableauside}\kern-\tableaurule}}}
\def\squares#1{\hbox{\count0=#1\noindent\loop\sqr
  \advance\count0 by-1 \ifnum\count0>0\repeat}}
\def\tableau#1{\vcenter{\offinterlineskip
  \tableaustep=\tableauside\advance\tableaustep by-\tableaurule
  \kern\normallineskip\hbox
    {\kern\normallineskip\vbox
      {\gettableau#1 0 }%
     \kern\normallineskip\kern\tableaurule}%
  \kern\normallineskip\kern\tableaurule}}
\def\gettableau#1 {\ifnum#1=0\let\next=\null\else
  \squares{#1}\let\next=\gettableau\fi\next}
\begin{document}

\title{Local curves, wild character varieties, and degenerations}
\author{
Duiliu-Emanuel 
Diaconescu,\\
 NHETC, Rutgers University,\\
 126 Frelinghuysen Road, Piscataway NJ 08854, USA\\
}
\date{}
\maketitle

\begin{abstract} 
Conjectural results for cohomological invariants of wild character varieties are obtained by counting curves in 
degenerate Calabi-Yau threefolds. A conjectural formula for $E$-polynomials 
is derived from the Gromov-Witten theory of local Calabi-Yau threefolds with normal crossing singularities. 
A refinement is also conjectured, generalizing 
existing results of Hausel, Mereb and Wong as well as recent joint work of Donagi, Pantev and the author for 
weighted Poincar\'e polynomials of wild character varieties.  
\end{abstract}

\tableofcontents

\section{Introduction} 

The cohomology of character varieties and Higgs bundle moduli spaces
is an important problem in geometry, topology and mathematical 
physics. Several approaches have been developed so far employing different methods. Arithmetic methods for character varieties have been 
used in \cite{HRV,HLRV,Arithmetic_wild, Zariski_closures} while similar methods for Higgs bundles
have been employed in  \cite{Ind_vb_higgs,Counting_Higgs,Ind_obj,Higgs_ind_proj_line}. 
The motive of the moduli space of irregular connections as well as irregular Higgs bundles over arbitrary fields has been determined in  \cite{Mot_conn_higgs}. Moreover, a  different approach based on wallcrossing for
linear chains  has been used in \cite{Mot_chains,y_genus_higgs} for motivic and Hirzebruch genus computations. At the same the topology of wild character varieties has been related to ploynomial invariants of legendrian knots in \cite{Fukaya_knots,Cluster_legendrian}. 
Finally, a string theoretic 
framework for this problem has been developed in \cite{wallpairs,BPSPW,Par_ref} based on an identification 
of perverse Betti numbers of Higgs moduli spaces with degeneracies 
of spinning BPS states in M-theory. In particular the conjectural formulas derived by Hausel and Rodriquez-Villegas \cite{HRV}, Hausel, Letellier and 
Rogriguez-Villegas \cite{HLRV} are identified with Gopakumar-Vafa expansions in the refined stable pair theory of certain Calabi-Yau threefolds. 
Furthermore, this framework also provides important evidence for the 
$P=W$ conjecture formulated by de Cataldo, Hausel and Migliorini 
\cite{Hodge_maps}. 

The string theoretic approach has been recently generalized to wild character varieties in \cite{BPS_wild}, providing a physical derivation and a generalization for the results of Hausel, Mereb and Wong \cite{Arithmetic_wild}. While the geometric framework 
and the spectral construction of \cite{BPS_wild} are general, explicit formulas are obtained only for wild character varieties with one singular point on the projective line. The goal of the present paper is to further extend these results to wild character varieties with multiple singular points on higher genus curves.

\subsection{Wild character varieties}\label{wildchar}
Wild character varieties are moduli spaces of Stokes data of irregular connections on curves 
with fixed irregular type. Such moduli spaces played a central role in Witten's work on geometric Langlands correspondence \cite{GL_wild} and were rigorously constructed by  Boalch in \cite{Quasi_hamiltonian, Braiding_stokes}.
The setup consists of  a smooth projective curve $C$ with $m$ pairwise distinct marked points $p_a$, $1\leq a \leq m$, and a complex reductive algebraic group $G$, which in this paper will be the general linear group, $GL(r, \IC)$, $r\geq 1$. 

Let 
${\bf T}_r\subset GL(r, \IC)$ be the standard maximal torus, and 
${\bf t}_r \subset {\mathfrak gl}(r,\IC)$ the corresponding Cartan subalgebra. 
As in \cite{Braiding_stokes}, a ${\bf t}_r$-valued irregular type 
at a point $p\in C$ is an equivalence class
\[ 
\sfQ \in {\bf t}_r({\widehat \CK}_p)/{\bf t}_r({\widehat \CO}_p),
\]
where ${\widehat \CO}_p$ is the natural completion of the local ring at $p$ and 
${\widehat \CK}_p$ is its field of fractions. 
Given a local coordinate $z$ centered at $p$, any irregular type 
$\sfQ$ admits a representative of the form 
\[ 
Q = \sum_{i=1}^{n-1} {A_i\over z^i} 
\]
for some integer $n\geq 2$ and some diagonal matrices $A_1, \ldots, A_{n-1}\in {\bf t}_r$. The common cetralizer of these Cartan elements depends only on the equivalence class 
$\sfQ$ and will be denoted by $H_\sfQ$. Since $A_1, \ldots, A_{n-1}$ are diagonal, $H_\sfQ$ is equivalent under conjugation with a canonical subgroup
$\times_{i=1}^\ell GL(m_i, \IC)\subset GL(r,\IC)$ for an ordered partition 
$r=m_1+\cdots + m_\ell$.  

Following \cite{Quasi_hamiltonian,Braiding_stokes}, a wild
curve is determined by the data $(C,p_a, \sfQ_a)$, $1\leq a\leq m$, 
where $\sfQ_a$ is a ${\bf t}_r$-valued irregular type at $p_a$.
For concreteness let $z_a$ be an affine coordinate at each $p_a$, and let 
\[
Q_a = \sum_{i=1}^{n_a-1} {A_{a,i}\over z_a^i} 
\]
be a representative of $\sfQ_a$ for each $1\leq a \leq m$.
Let $\sfQ = (\sfQ_1, \ldots \sfQ_m)$. 
The wild character variety associated to a wild curve is the moduli space of Stokes data for irregular flat connections on $C\setminus\{p_1, \ldots, p_m\}$ 
which are locally gauge equivalent to 
\[
dQ_a + {\rm terms\ of\ order}\ \geq -1 
\] 
in the infinitesimal neighhborhood of each marked point $p_a$. 
One of the main results of \cite{Braiding_stokes}
proves that this is a smooth quasi-projective variety with a natural Poisson structure. Moreover, the symplectic leaves 
of the Poisson structure are obtained by fixing the conjugacy class 
$\tau_a$ 
in $H_{\sfQ_a}$ of the formal monodromy of the irregular flat connections at 
each $p_a$ for all $1\leq a \leq m$. Note that the conjugacy class of the formal mondromy in 
$H_{\sfQ_a}$ coincides with the conjugacy class 
of the local monodromy around the puncture, which is determined by the residue of the irregular connection. 

Assuming the marked curve fixed, let $\CS_{\sfQ, \tau}$ be the symplectic leaf determined 
by $\tau=(\tau_1, \ldots, \tau_m)$.
Throughout this paper all irregular types will be 
 chosen such that 
\begin{itemize} 
\item[$(W.1)$]
$H_\sfQ=\times_{i=1}^\ell GL(m_i, \IC)\subset GL(r,\IC)$ and $H_\sfQ$ coincides with the centralizer of the 
leading term in the Laurent expansion. 
\item[$(W.2)$]
Each $\tau_a$, $1\leq a\leq m$,  is the conjugacy class of a diagonal matrix $\sfT_a$ with eigenvalues 
\be\label{eq:formalmon}
(\underbrace{\tau_{a,1}, \ldots, \tau_{a,1}}_{m_{a,1}}, \ldots, \underbrace{\tau_{a,\ell_a}, \ldots, \tau_{a,\ell_a}}_{m_{a,\ell_a}}),
\ee
where $\tau_{a,1}, \ldots, \tau_{a,\ell_a}$ are pairwise distinct complex numbers.
\end{itemize}
 For each $a$ let $\mu_a$ denote the partition of $r$ determined by $(m_{a,1}, \ldots, m_{a, \ell_a})$ and let $\umu=(\mu_1, \ldots, \mu_m)$. 

As shown in \cite[Thm. 8.2]{Braiding_stokes}, for fixed irregular types 
$\sfQ=(\sfQ_1, \ldots, \sfQ_m)$ and fixed conjugacy classes 
$\tau=(\tau_1,\ldots, \tau_m)$ the associated moduli space $\CS_{\sfQ,\tau}$ of 
Stokes data has an explicit presentation as an affine algebraic 
quotient. For sufficiently generic $\tau$ this moduli space is a smooth quasi-projective variety equipped with a holomorphic symplectic structure. Its dimension is given by 
\be\label{eq:dformula}
d(\umu, \un, g)= 2r^2(g-1)+\sum_{a=1}^m n_a(r^2-\sum_{i=1}^{\ell_a} \mu_{a,i}^2) +2.
\ee
where $\un= (n_1,\ldots, n_m)$.

\subsection{Higgs bundles and wild non-abelian Hodge correspondence}\label{irreghiggs}
Results proven by Biquard and Boalch \cite{Wild_curves} and Sabbah \cite{Harmonic_metrics} show that
moduli spaces of irregular filtered flat connections are related by hyper-K\"ahler 
rotations to moduli spaces of irregular parabolic Higgs bundles. 
This relation is commonly refered to as the wild non-abelian Hodge correspondence. A very clear and explicit statement can be found in 
\cite{HK_wild} to which the reader is refered for more details. 

For the purposes of the present paper, it suffices to note that given sufficiently generic parameters $(\sfQ,\tau)$ the resulting Higgs bundle moduli problems 
 admit a specific formulation leading to a natural connection to Calabi-Yau enumerative geometry. This construction as well as its relation to 
the moduli spaces of \cite{HK_wild} is explained in detail in Sections 
2.1 and 2.3 of \cite{BPS_wild}. 
Very briefly, one employs Higgs bundles with a coefficient line bundle 
$M= K_C(D)$ where $D=\sum_{a=1}^m n_a p_a$ is the total polar divisor of the 
fixed irregular types.  The parabolic structure consists of  a locally-free filtration along each non-reduced divisor $D_a =n_ap_a$, $1\leq a\leq m$. 
The numerical type of such a filtration is specified by 
a sequence of positive integers $(m_{a,1}, \ldots, m_{a, \ell_a})$, 
where $m_{a,i}$ is the length of the $i$-th successive quotient 
of the $a$-th filtration  as an $\CO_{D_a}$-module. Let 
$\um = (m_{a,i})$, $1\leq i \leq \ell_a$, $1\leq a \leq m$.
The Higgs field is required to preserve all these filtrations, 
and moreover the associated graded ${\rm gr}\big(\Phi|_{D_a})$ is required to act as a scalar on each succesive quotient. 
More precisely, for each $1\leq a\leq m$ one requires 
\[ 
{\rm gr}\big(\Phi|_{D_a}) = \oplus_{i=1}^{\ell_a} \xi_{a,i} \otimes 
{\bf 1}_i
\]
where 
$\xi_{a,i}\in H^0(D_a, M|_{D_a})$, $1\leq i \leq \ell_a$ are some fixed sections of $M$ over $D_a=n_ap_a$. Let $\uxi_a=(\xi_{a,i})$ with $1\leq i \leq \ell_a$ and let also $\uxi=(\xi_{a,i})$, $1\leq i \leq \ell_a$, $1\leq a \leq m$. Throughout this paper the sections $\uxi_a$ will be 
assumed generic for all $1\leq a\leq m$, meaning that their values 
$\xi_{a,i}(p_a)$ at the reduced marked point are pairwise distinct and nonzero. 

In addition one has to specify real parabolic weights $\ualpha = (\alpha_{a,i})$, $1\leq i \leq \ell_a$, $1\leq a \leq m$ and impose the standard parabolic stability conditions on such objects as in \cite{modpar}.
This results in a moduli stack 
${\mathfrak H}^{ss}_\uxi(C,D; \ualpha, \um, d)$ of semistable objects, where $d$ is the degree of the Higgs bundles. Moreover, strictly semistable objects are absent for
sufficiently generic parabolic weights $\ualpha$. In such cases the moduli stack is a $\IC^\times$-gerbe over a smooth quasi-projective 
moduli space $\CH^{s}_\uxi(C,D; \ualpha, \um, d)$ as usual. 

To conclude, note that the wild non-abelian Hodge correspondence leads to a wild variant of the $P=W$ conjecture of de Cataldo, Hausel and Migliorini \cite{Hodge_maps}. This correspondence 
identifies the weight filtration on the cohomology of 
a character variety and the perverse Leray filtration on the cohomology of 
the associated Higgs bundle moduli space.  

\subsection{Spectral correspondence and 
Calabi-Yau threefolds}\label{spcorresp}

As shown in detail in \cite[Sect. 3.1]{BPS_wild}, the string theoretic construction
is based on a 
spectral correspondence relating the irregular Higgs bundles to 
dimension one
sheaves on a holomorphic symplectic surface $S_\uxi$. Employing
the construction of \cite{structures}, this surface is obtained by blowing up the images of the sections $\xi_{a,i}$, taking a minimal resolution of singularities, 
and removing an anti-canonical divisor. 

The compact curve classes on $S_\uxi$ are in one-to-one correspondence with collections
non-negative integers $\um=(m_{a,i})$, $1\leq i\leq \ell_a$, 
$1\leq a\leq m$ such that the sum $\sum_{i=1}^{\ell_a} m_{a,i}$
takes the same value, $r\geq 1$, for each $1\leq a\leq m$. 
Any compact curve in such a class is a finite $r:1$ cover of $C$. The curve class corresponding to $\um=(m_{a,i})$ will be denoted by $\beta(\um)$. 

As shown in 
\cite[Sect.3.2-3.4]{BPS_wild}, for sufficiently generic $\uxi$ 
 the moduli stack ${\mathfrak H}^{ss}_\uxi(C,D; \ualpha, \um, d)$  is isomorphic to a moduli stack of Bridgeland semistable pure dimension one sheaves $F$ on $S_\uxi$  with compact support.  The topological invariants of these sheaves are 
\[ 
\ch_1(F) = \beta(\um), \qquad \chi(F) =  d - r(g-1)
\]
while the Bridgeland stability condition is determined by the parabolic weights $\ualpha$. 

This correspondence can be lifted one step further. Namely 
let $Y_\uxi$ be the total space of the canonical bundle $K_{S_\uxi}$, which is isomorphic to $S_\uxi\times \IA^1$. 
Then any moduli stack of compactly supported Bridgeland stable 
pure dimension one sheaves on $Y_\uxi$ is isomorphic to a product of the form ${\mathfrak H}^{s}_\uxi(C,D; \ualpha, \um, d)\times \IA^1$. 
The threefold $Y_\uxi$ will be called a wild local curve throughout this paper. Occasionally it  will be also referred to as a spectral threefold for ease of exposition. 

Applying the mathematical framework of \cite{HST,GV_vanishing}, in the present context the degeneracies of spinning BPS particles in M-theory are given by  the perverse Betti numbers of the Higgs bundle moduli spaces ${\mathfrak H}^{s}_\uxi(C,D; \ualpha, \um, d)$. 
Therefore, using the $P=W$ conjecture and the 
 refined Gopakumar-Vafa expansion \cite{GV_II,spinBH,Mtop,ref_vert}, 
it follows that the weighted Betti numbers of wild character varieties
are completely determined by the refined 
stable pair theory of $Y_\uxi$. 

\subsection{Curve counting invariants and degenerations}\label{counting} 

Computing the refined stable pair theory of the threefolds $Y_\uxi$ 
constructed above is a difficult and challenging problem. 
Conjectural formulas were derived in \cite{BPS_wild} for genus zero curves with a single marked point. The derivation in loc. cit. uses refined virtual localization as in \cite{Membranes_sheaves} and the refined colored 
variant \cite{homflypairs} of the Oblomkov-Shende-Rasmussen conjectures \cite{OS,ORS} to reduce the problem to refined torus link invariants. The latter are computed in 
turn using refined Chern-Simons theory \cite{refCS} and large $N$ 
duality. 

This approach cannot be implemented directly for a genus $g$ 
curve with multiple marked points since there is no torus action on the associated spectral surface $S_\uxi$. In fact, given the absence of a torus action, even the computation of unrefined invariants poses significant problems. 
The present work develops an alternative strategy employing 
nodal degenerations of the curve $C$ to construct a normal crossing 
spectral surface equipped with a suitable torus action. This torus action 
is then used to find an explicit formula for the unrefined stable pair theory 
of the resulting degenerate threefolds, providing conjectural formulas for 
 $E$-polynomials of wild character varieties. The resulting formulas admit a further conjectural refinement which leads to explicit conjectures for 
weighted Poincar\'e polynomials. In more detail, the main steps of this approach  are summarized below. 

$(i)$ As defined by Pandharipande and Thomas \cite{stabpairsI}, a stable pair
$(F,s)$ on the smooth spectral threefold $Y_\uxi$ is a pure dimension one sheaf $F$ equipped with a generically surjective section $s:\CO_{Y_\uxi}\to F$. 
The support of $F$ is required to be compact and one fixes the invariants $\ch_1(F)=\beta(\um)\in H_2(Y_\uxi,\IZ)$
and $c = \chi(F)$. The resulting moduli space of pairs ${\mathcal P}(Y_\uxi, \um,c)$ has a perfect obstruction theory and a virtual cycle of degree zero. However, since it is non-compact and there is no torus action with compact fixed locus, hence stable pair invariants cannot be defined by virtual 
integration. One has to employ instead Behrend's constructible function approach  \cite{micro}, definining the  unrefined invariants as weighted 
Euler characteristics of moduli spaces, 
\be\label{eq:stpairA}
PT(Y_\uxi, \um, c) = \chi( {\mathcal P}(Y_\uxi, \um,c), \nu^B).
\ee
The weights are encoded in the constructible function  $\nu^B: {\mathcal P}(Y_\uxi, \um,c)\to \IZ$ constructed in \cite{micro}. 
As opposed to virtual integration, this definition is not known to be deformation invariant in general. However, in the present context, it will be assumed  
that the constructible invariants are in fact invariant 
under deformations of $Y_\uxi$ induced by deformations of 
the base curve $C$ and the marked points. 

$(ii)$  A nodal degeneration $\oC$ of the curve $C$ is constructed in 
Section \ref{degsect} consisting of a smooth central component $C_0$ 
and $m$ pairwise disjoint rational components $C_1, \ldots, C_m$. 
Each projective line $C_a$ contains exactly one marked point $p_a$ and 
intersects $C_0$ transversely at another point $\nu_a$. For fixed 
data $\uxi$, this degeneration yields  normal crossing degenerations of the spectral surface $S_\uxi$ as well as the threefold $Y_\uxi$. 
Moreover, as shown in Section \ref{toract}, for a suitable choice of data, the resulting normal crossing threefold admits a torus action preserving the normal crossing Calabi-Yau structure. 
Using the work of J. Li and B. Wu \cite{Good_deg},  for any pair $(\um, c)$ there is a moduli
space of stable pairs on the normal crossing threefold $\oY_\uxi$ equipped with a 
perfect obstruction theory.
Moreover, the fixed locus of the torus action is compact, hence 
in this limit one can construct equivariant residual stable pair invariants
$PT(\oY_\uxi, \um, c)$. Then a stronger deformation invariance  assumption will be made in this paper, claiming that 
\be\label{eq:stpairB}
PT(\oY_\uxi,\um, c) = PT(Y_\uxi, \um,c)
\ee
for all $\um,c$. 

$(iii)$ The equivariant residual stable pair theory of degenerations is explicitely computed in Section \ref{GWwildlocal} using GW/PT correspondence, which will 
be assumed without proof in the present context. 
In principle one should be able to give a proof 
by analogy with \cite{Curves_modular}, but this would be beyond the scope of the present paper.  Granting this correspondence, the residual Gromov-Witten theory of degenerations is computed in Section \ref{GWwildlocal} 
using J. Li's theory of relative stable maps 
\cite{Stable_sing} and the degeneration formula \cite{Deg_formula}, 
the relative virtual localization theorem of Graber and Vakil \cite{Relative_loc},  as well as the local TQFT formalism of Bryan and Pandharipande \cite{Loc_curves}. 
The computation also uses the Marino-Vafa 
formula proven by C.-C. M. Liu, K. Liu
and J. Zhou \cite{MV_formula} respectively Okounkov and Pandharipande \cite{Hodge_unknot}
as well as results of Bryan and Pandharipande 
\cite{Curves_TQFT}, Okounkov and Pandharipande 
\cite{GW_completed,Virasoro_curves} and 
C.-C. M. Liu, K. Liu
and J. Zhou \cite{Two_partition}
 on relative  Gromov-Witten invariants of (local) curves. For completeness, a 
self-contained exposition of relative stable maps and relative virtual 
localization, including one the main examples used in the computation, is 
provided in Section \ref{relativemaps}. 
The final outcome is presented below. 

\subsection{The main formula} 
To summarize the current setup, in this section ${\oC}$ will be a nodal curve consisting of a smooth 
genus $g$ central component $C_0$ and $m$ pairwise disjoint 
smooth genus zero components $C_1, \ldots, C_m$.  Each component 
$C_a$ intersects $C_0$ transversely at one point $\nu_a$ and. One also 
specifies a nonreduced effective divisor $D=\sum_{a=1}^m n_ap_a$ on 
$\oC$ such that $p_a\in C_a\setminus \{\nu_a\}$ for each 
$1\leq a\leq m$. Each marked point is assigned in addition 
a collection $\uxi_a= (\uxi_{a,1}, \ldots, \uxi_{a,\ell_a})$ of $\ell_a\geq 1$ 
pairwise disjoint, nonzero, sections of $M=\omega_{\oC}(D)$ over the thickening $D_a=n_ap_a$, where $\omega_{\oC}$ is the dualizing sheaf of $\oC$.
 
Using this data one constructs a normall crossing spectral surface $\oS_\uxi$ following the algorithm of 
\cite{structures}, as in \cite[Sect. 3.1]{BPS_wild}. One has to blow up the images of the sections $\uxi$ in the total 
space of $M$, take a minimal resolution of singularities, and remove a certain divisor. The spectral threefold $\oY_\uxi$ is the total space of the dualizing sheaf of $\oS_\uxi$, hence it has a normal-crossing Calabi-Yau structure. It consists of a central component, $Y_0$, which is isomorphic to the total space of a rank two bundle over $C_0$ and $m$ components 
$Y_a$, $1\leq a\leq m$, each of them isomorphic to the total space of a line bundle over a surface $S_a$, $1\leq a\leq m$. By construction,  both $Y_a$ and $S_a$ have natural log Calabi-Yau structures for all $1\leq a\leq m$.  The polar divisor of the dualizing sheaf of $Y_a$, denoted by
$\Delta_a$, $1\leq a\leq m$, is isomorphic to the affine plane. Moreover, 
the threefolds $Y_a$ are glued to the central component $Y_0$ along the 
divisors $\Delta_a$, which are naturally identified with fibers of the projection 
$Y_0\to C_0$.  Each pair $(Y_a,\Delta_a)$ will be called a wild cap 
in the following. 

As shown in Section \ref{toract}, assuming that 
\[ 
n_1=\cdots = n_m =n,
\]
and making a suitable choice of local data $\uxi$, one constructs a torus 
action on $\oS_\uxi$.  This action lifts to an action on the 
spectral threefold $\oY_\uxi$ preserving the 
normal-crossing Calabi-Yau structure. 
Moreover for each $1\leq a\leq m$ the cone of compact curve classes 
on $Y_a$ is freely generated by $\ell_a$ torus invariant 
sections of the projection map  $Y_a\to C_a$.  Therefore the compact curve classes 
on $Y_a$ are in one-to-one correspondence with collections 
$(m_{i,a})\neq (0,\ldots,0)$, 
$1\leq i \leq \ell_a$ of non-negative integers. At the same time the cone of  compact curve classes in $Y_0$ is freely generated by the zero section.
Hence any such class is classified by the degree $r\geq 1$.  In conclusion 
the compact curve classes on $\oY_\uxi$ in one-to-one correspondence with  numerical data $(\um, r)$. 

As shown in \cite{MNOP_I,MNOP_II,stabpairsI}, in the context of
GW/DT/PT correspondence one has to consider stable maps to $\oY_\uxi$ with disconnected domains such that no connected component is mapped 
to a point.
Using the degeneration formula proven in \cite{Deg_formula} it follows 
that  only pairs $(\um, r)$ satisfying 
\be\label{eq:sumcond} 
\sum_{i=1}^{\ell_a} m_{a,i} =r \geq 1
\ee
for all $1\leq a\leq m$ contribute to the Gromov-Witten theory of $\oY_\uxi$. In particular $r$ is determined by $\um$, hence the residual equivariant invariants will be denoted by $GW^\bullet_{\um,h}$, 
where $h\in \IZ$ is the arithmetic genus of the domain. By convention 
for any $h \in \IZ$, and any integral vector $\um$ 
let $GW^\bullet_{\um,h}=0$ unless $\um$ has non-negative entries satisfying condition \eqref{eq:sumcond}.  
Therefore the partition function 
takes the form 
\be\label{eq:GWfct}
Z(\oY_\uxi; {\sf x}, g_s) = 1+\sum_{\um}\ \sum_{h\in \IZ}\ 
g_s^{2h-2} GW_{\um,h} \prod_{a=1}^m \prod_{i=1}^{\ell_a} x_{a,i}^{m_{a,i}}.
\ee

An explicit formula for $Z(\oY_\uxi; {\sf x}, g_s)$ is derived in Section \ref{GWwildlocal} by relative virtual localization.
The formula is written as a sum over partitions 
\[
\lambda = (\lambda_1\geq \lambda_2 \geq \cdots\geq \lambda_{l(\lambda)}),
\]
 identified with Young diagrams consisting of $l(\lambda)$ left-aligned rows such that the $i$-th row contains $\lambda_i$ boxes.  The following notation and definitions are used in the computation. 
\begin{itemize}
\item The total number of boxes in $\lambda$ will be denoted by $|\lambda|$ and its transpose will be denoted by  $\lambda^t$.
\item The content of $\lambda$ is defined as 
\[ 
c(\lambda) = \sum_{\Box\in \lambda} (a(\Box)-l(\Box) ) 
\]
where $a(\Box), l(\Box)$ are the arm and leg length respectively.  
\item The Schur function corresponding to $\lambda$ will be denoted by $s_\lambda({\sf x})$. The fusion coefficients for $\ell$ partitions $\nu_1, \ldots, \nu_\ell$ are defined by 
\[ 
\prod_{i=1}^\ell s_{\nu_i}({\sf x}) = \sum_{|\lambda|=|\nu_1|+\cdots+ |\mu_\ell|} c^{\lambda}_{\nu_1, \ldots, \nu_\ell} s_\lambda({\sf x}). 
\]
\end{itemize} 
Then the formula derived in Section \ref{GWwildlocal} reads 
\be\label{eq:wildformulaA} 
\bal 
Z_{GW}(\oY_\uxi; {\sf x}, g_s)  = 1+\sum_{|\lambda|>0} s_{\lambda^t}(\uq)^{2-2g-m} \prod_{a=1}^m F_{n-1,\ell_a,\lambda}({\sf x}_a,q) 
\eal
\ee
where $q=e^{ig_s}$, $\uq=(q^{1/2},q^{3/2}, \ldots)$, and 
\be\label{eq:wildfactorsA}
F_{k,\ell,\lambda}({\sf x}, q) = q^{-kc(\lambda)}\sum_{|\nu_1|+\cdots+|\nu_\ell|=|\lambda|} 
c^\lambda_{\nu_1, \ldots, \nu_{\ell}} \prod_{i=1}^{\ell} x_{i}^{|\nu_i|}q^{kc(\nu_i)} s_{\nu_i^t}(\uq).
\ee
for any integers $k,\ell\geq 1$. The factors $F_{n-1,\ell_a, \lambda}({\sf x}_a, q)$ encode the relative Gromov-Witten invariants 
of the wild caps $(Y_a, \Delta_a)$, which are the main novelty in this 
computation. 

Assuming GW/PT correspondence, the above partition function coincides with the generating function of residual stable pairs invariants $PT(\oY_\uxi,\um,c)$ of the threefolds $Y_\uxi$. 
Further assuming deformation invariance as in\eqref{eq:stpairB}, the above formula provides an explicit expression for the 
generating function 
\[ 
Z_{PT}(Y_\uxi, {\sf x}, q) = 1+ \sum_{\um}\ \sum_{c\in \IZ} PT(Y_\uxi, \um, c) 
(-q)^c \prod_{a=1}^m \prod_{i=1}^{\ell_a} x_{a,i}^{m_{a,i}}
\]
associated to a smooth spectral threefold $Y_\uxi$. 

The above formula is derived under the assumption that all integers $n_a$ used as input data
in the construction of $Y_\uxi$ 
are equal,
$n_1=\cdots = n_m$. Moreover, one also requires the sections $\uxi_a$ 
to satisfy an equivariance condition for all $1\leq a \leq m$ in order 
 for a torus action to exist in the degenerate limit. 
At the same time, the invariants \eqref{eq:stpairA} are defined more generally for any values of $n_a$ and for generic sections $\uxi_a$, 
$1\leq a\leq m$. 
The next section formulates a conjectural refined generalization 
of formula \eqref{eq:wildformulaA} 
 allowing arbitrary  values of $n_a$, as well as arbitrary generic sections 
$\uxi_a$, 
$1\leq a\leq m$.

\subsection{Conjectural refinement}
The conjectural refinement of the partition function \eqref{eq:wildformulaA}
is inferred from the structure of the wild factors \eqref{eq:wildfactorsA}
by comparison with previous refined formulas derived in  \cite{BPSPW,Par_ref,BPS_wild}. 
The refinement will be written in terms of Macdonald polynomials $P_\lambda(s,t;{\sf x})$ and uses the following quantities:
\begin{itemize} 
\item the fusion coefficients $N^\lambda_{\nu_1, \ldots, \nu_\ell}(s,t)$ defined by the product identity 
\[
\prod_{i=1}^\ell P_{\nu_i}(s,t; {\sf x}) = \sum_{|\lambda|=|\nu_1|+\cdots + |\nu_\ell|} N^\lambda_{\nu_1, \ldots, \nu_\ell}(s,t) P_\lambda(s,t; {\sf x}) 
\]
\item the specializations 
\[
R_\mu(q,y) = P_\mu(s,t; \ut)\big|_{s=qy,\ t=qy^{-1}}, \qquad
L_\mu(q,y) = P_\mu(t,s;\us)|_{s=qy,\ t=qy^{-1}}
\]
\[ 
{\widetilde N}^\lambda_{\nu_1, \ldots, \nu_\ell}(q,y) = 
N^\lambda_{\nu_1, \ldots, \nu_\ell}(qy,qy^{-1})
\]
where 
\[
\ut = (t^{1/2}, t^{3/2}, \ldots), \qquad \us= (s^{1/2}, s^{3/2}, \ldots),
\]
and
\item the framing factors 
\[
g_\mu(q,y) = f_\mu(qy,qy^{-1}), \qquad f_\mu(s,t) = \prod_{\Box\in \mu} s^{a(\Box)} t^{-l(\Box)}.
\]
\end{itemize} 
One of the main building blocks of the refined formulas derived in 
\cite{BPSPW,Par_ref} is the polynomial function 
\[
T_{g,\lambda}(q,y) = \prod_{\Box\in \lambda} 
(qy)^{-(2a(\Box)+1)g}(1-y^{a(\Box)-l(\Box)} q^{h(\Box)})^{2g}
\]
encoding the genus dependence of the partition function. 
Furthermore, by comparison with the refined formula of \cite{BPS_wild}, 
one is led to conjecture the following refinement of the wild factor 
\eqref{eq:wildfactorsA}. 
\be\label{eq:wildcapD} 
\bal 
{G}_{k,\ell, \lambda}({\sf x}, q,y) = \sum_{\substack{ 
\mu_1, \ldots,\mu_\ell\\ |\mu_1|+\cdots + |\mu_\ell|=|\lambda|\\}} {\widetilde N}^\lambda_{\mu_1, \ldots, \mu_\ell}(q,y) g_\lambda(q,y)^{-k} \prod_{i=1}^\ell \left(x_i^{|\mu_i|} g_{\mu_i}^{k}(q,y) L_{\mu_i^t}(q,y) \right).
\eal 
\ee
Then the refined stable pair theory of a threefold $Y_\uxi$ with arbitrary data $(n_a, \ell_a, \uxi_a)$, $1\leq a \leq m$, at the marked points is  conjectured to be 
\be\label{eq:wildformulaB}
\bal
Z_{PT}^{ref}(Y_\uxi; & q,y, {\sf x}) = \\
& 1+\sum_{|\lambda|>0} 
g_\lambda(q,y)R_\lambda(q,y) T_{g,\lambda}(q,y)
L_{\lambda^t}(q,y)^{1-m} \prod_{a=1}^m G_{n_a-1,\ell_a,\lambda}({\sf x}_a, q,y)\\
\eal 
\ee
As a consistency check, this formula reduces to the one conjectured in
\cite{BPSPW} in the absence of marked points, in which case the wild cap factors are absent. Moreover, it also reduces to the main formula conjectured in \cite{BPS_wild} for a genus zero curve with one marked point. 
Moreover hard numerical evidence for this conjecture is obtained by comparison with the 
main conjecture of Hausel, Mereb and Wong \cite{Arithmetic_wild}, as discussed 
in the next section.

\subsection{Conjectures for wild character varieties}\label{wildconjectures}
As explained in Section \ref{irreghiggs}, for sufficiently generic data $\uxi$ the wild non-abelian Hodge correspondence relates moduli spaces of irregular parabolic Higgs bundles with singular type $\uxi$ to moduli spaces of filtered flat connections. A set of generic of local data $\uxi$ and a set of generic parabolic 
weights determine uniquely a collection $\Gamma(\uxi)$ of irregular types as well 
as a collection $\tau$ of conjugacy classes 
associated to the marked points. Both $\Gamma(\uxi)$ and $\tau$ satisfy 
conditions $(W.1)$ and $(W.2)$ in Section \ref{wildchar}. 
This relation leads to a wild variant of the 
$P=W$ conjecture of \cite{Hodge_maps} identifying the perverse Leray filtration on the moduli space of Higgs bundles with the weight filtration 
on the cohomology of wild character varieties. 
Using the refined Gopakumar-Vafa expansion, the refined stable pair 
partition function \eqref{eq:wildformulaB} 
yields explicit predictions for weighted Poincar\'e polynomials of wild character varieties. 
By analogy with \cite{BPSPW,Par_ref,BPS_wild}, in the present context the 
refined Gopakumar-Vafa formula is a conjectural expansion 
\be\label{eq:refGVa}
\bal
 {\rm ln}\, Z_{ref}({Y_\uxi}; & {\sf x},q,y) = \\
&  -
\sum_{k\geq 1}\ \sum_{|\mu_1|=\cdots =|\mu_m|=r\geq 1} 
\ \prod_{a=1}^m
{m_{\mu_a}({ x}^k_a)
\over k}
 {y^{-kr}
(qy^{-1})^{kd(\umu,\un,g)/2} P_{\umu,\un}((qy)^{-k}, y^k)) \over 
 (1-(qy)^{-k})(1-(qy^{-1})^{k})}\\
\eal
\ee
where $P_{\umu,\un}(u,v)$ are polynomials with integer coefficients 
which count the degeneracies of spinning BPS states in M-theory. 
 The coefficients $m_\mu({\sf x}^k)$ are the monomial symmetric functions
evaluated at ${\sf x}^k=(x_1^k, x_2^k, \ldots)$ for any $k\geq 1$. Note that ${\sf x}_a=(x_1, \ldots, x_{\ell_a},0,0, \ldots)$ for each 
$1\leq a\leq m$.

Using the mathematical theory of \cite{HST,GV_vanishing}, since $Y_\uxi$ is a spectral threefold,  the polynomials 
$P_{\umu,\un}$ are identified with the perverse Poincar\'e 
polynomials of the Higgs bundle moduli spaces $\CH^{s}_\uxi(C,D; \ualpha, \um, d)$ for sufficiently generic $\ualpha$. 
This yields the following conjectural statement for wild character varieties. 

Let $\CS_{\Gamma(\uxi), \tau}$ be the wild character variety corresponding to 
some generic sections $\uxi$ and generic conjugacy classes $\tau$ 
satisfying conditions $(W.1)$, $(W.2)$ in Section \ref{wildchar}. Let 
\[ 
P(\CS_{\Gamma(\uxi), \tau}; u,v) = \sum_{i,j} {\rm dim}\, {\rm Gr}_i^W
H^j(\CS_{\Gamma(\uxi), \tau}) u^{i/2} (-v)^j
\]
be the weighted Poincar\'e polynomial of $\CS_{\Gamma(\uxi), \tau}$ 
where ${\rm Gr}_i^W
H^j(\CS_{\Gamma(\uxi), \tau})$ are the successive quotients of the weight filtration on cohomology. Then 
\[ 
P(\CS_{\Gamma(\uxi), \tau}; u,v)  = P_{\umu, \un}(u,v). 
\]
In particular, under the current genericity assumptions the weighted Poincar\'e polynomial depends only on the discrete data $(\umu,\un)$. 

The main supporting evidence for this conjecture is obtained from the comparison with the conjecture of Hausel, Mereb and Wong \cite{Arithmetic_wild} which yields explicit predictions for 
partitions of the form $\umu= ((1^r), \ldots, (1^r))$.  
Their partition function is defined as  
\[ 
Z_{HMW}(z,w) = 1+\sum_{|\lambda|>0} \Omega_{\lambda}^{g,n}(z,w) \prod_{a=1}^m
{\wH}_\lambda({\sf x}_a; z^2, w^2) 
\]
where:
\begin{itemize}
\item the sum in the right hand side is over all Young diagrams 
$\lambda$ with $|\lambda|>0$, 
\item 
for each $\lambda$ 
\[ 
\Omega_\lambda^{g,n} (z,w) = \prod_{\Box \in \lambda} 
{ (-z^{2a(\Box)}w^{2l(\Box)})^{n-m}
(z^{2a(\Box)+1} - w^{2l(\Box)+1})^{2g} \over 
 (z^{2a(\Box)+2} -w^{2l(\Box)})(z^{2a(\Box)} -w^{2l(\Box)+2})},
\] 
\item ${\sf x}_a = (x_{a,1}, x_{a,2}, \ldots)$ is an infinite set of formal variables for each $1\leq a\leq m$, and ${\wH}_\lambda({\sf x}_a; z^2,w^2)$ are the 
modified Macdonald polynomials. 
\end{itemize}
Let ${\mathbb H}_{\mu,n}(z,w)$ be defined by 
\be\label{eq:HMWa}
{\rm ln}\, Z_{HMW}(z,w) = \sum_{k\geq 1}\sum_{\mu} {(-1)^{(n-m)|\mu|}
w^{kd(\mu,n,g)}\, {\mathbb H}_{\umu,n}(z^k,w^k)\over (1-z^{2k})(w^{2k}-1) }\prod_{a=1}^m 
m_{\mu_a}({x}_a^k) 
\ee
where the sum is again over all Young diagrams, $m_{\mu}({\sf x})$ are the monomial symmetric functions and ${\sf x}^k = (x_1^k, x_2^k, \ldots)$.  
Then for $\mu_a=(1^r)$, $1\leq a\leq r$, one has the following conjectural formula 
\be\label{eq:HMWb}
WP(\CS_{{\sfQ},\sfT}; u,v) = {\mathbb H}_{(1^r),n}
(u^{1/2}, u^{-1/2}v^{-1})
\ee
where $n= \sum_{a=1}^m n_a$. 
Note that the $v=1$ specialization of this conjecture is proven in \cite{Arithmetic_wild} using arithmetic methods. 

Although the two conjectural formulas are quite different, direct numerical comparison shows that they yield the same predictions for 
\begin{itemize} 
\item $0\leq g \leq 2$, $m=2$, $\mu_1=\mu_2=(1^3)$, $2\leq {\rm deg}(D)=n_1+n_2\leq 7$
\item $0\leq  g \leq 1$, $m=3$, $\mu_1=\mu_2=\mu_3=(1^3)$, $3\leq {\rm deg}(D)=n_1+n_2+n_3\leq 9$. 
\end{itemize} 
It should be noted that these are highly nontrivial tests. For example 
for $g=2$, ${\deg(D)}=7$, $\mu_1=\mu_2=(1^3)$, as well as  
$g=1$, ${\deg(D)}=9$, $\mu_1=\mu_2=\mu_3=(1^3)$ 
the polynomials $P_{\umu, n}(u,v)$ 
have bidegree $(56,56)$, hence a total of 3248 terms. 
Moreover, as a further consistency check, explicit computations for 
rank $r=3$ and $r=4$ character varieties confirm that the partition 
function \eqref{eq:wildformulaB} has indeed a BPS expansion of the form 
\eqref{eq:refGVa} for arbitrary partitions $\mu_a$, including   $\mu_a \neq (1^r)$.  For illustration some numerical results for the polynomials $P_{\umu, \un}$ are listed in Appendix \ref{examples}. 
\bigskip

{\bf Acknowledgments.} I am very grateful to Ron Donagi and Tony Pantev for collaboration on the related project \cite{BPS_wild} and to 
Philip Boalch, Tamas Hausel, Davesh Maulik and Yan Soibelman for very helpful discussions and comments on the manuscript. During the completion of this work the author was partially supported by NSF grant DMS-1501612.

\section{Wild degenerate local curves}\label{degwild}

\subsection{Degenerations}\label{degsect}
In this section $\oC$ will be a nodal curve 
consisting of $m$ rational components 
$C_1, \ldots, C_{m}$ and a smooth projective genus $g$ curve $C_0$. Each rational component $C_a$, $1\leq a\leq m$, intersects $C_0$ transversely at a node $\nu_a\in C_0$ and is at the same time disjoint from all other rational components $C_b$, $b\neq a$.  Moreover, each component $C_a$ contains a marked point $p_a\neq \nu_a$, and each marked point
is assigned the pair of positive integers $(n_a, \ell_a)$. This data will be denoted by 
$\un=(n_a)$, $\uell=(\ell_a)$ with $1\leq a\leq m$.
Furthermore let  $D_a=n_ap_a$, $1\leq a\leq m$, and $D=\sum_{a=1}^m D_a$.

A degenerate wild curve is then specified by the collection 
$(\oC, p_1, \ldots, p_m)$ and a collection $\uxi_a=(\xi_{a,i})$, $1\leq i \leq \ell_a$
of generic
sections of $M|_{D_a}$ for each $1\leq a \leq m$. 
Recall that the genericity condition requires the values $\xi_{a,i}(p_a)$ to be
pairwise distinct and nonzero for each $1\leq a \leq m$. Let $\uxi=(\uxi_{a,i})$, $1\leq i \leq \ell_a$, 
$1\leq a\leq m$. 

Let $M=\omega_{\oC}(D)$,
 where $\omega_{\oC}$ is the dualizing line bundle of $\oC$. 
In analogy with \cite{structures,BPS_wild}, one constructs a surface
${\oS}_\uxi$ by blowing-up
the total space of $M$ along the images of the sections $\uxi_{a,i}$, taking a minimal resolution of singularities, and removing an anticanonical divisor.
This process can be alternatively viewed as a sequence of 
successive blowups of M at closed points lying recursively on the strict transforms of the sections 
$\xi_{a,i}$. The resulting surface is reducible, with normal crossing singularities, and there is a natural projection 
map $\eta: \oS_\uxi \to \oC$. The smooth irreducible components of $\oS_\uxi$ are 
 $S_a=\eta^{-1}(C_{a})$ and $S_0=\eta^{-1}(C_0)$. The central component $S_0=\pi^{-1}(C_0)$ is isomorphic to the total space 
of $$M|_{C_0}\simeq K_{C_0}(\nu_1+\cdots +\nu_{m})$$
and intersects each $S_{a}$ transversely along the fiber $F_a=\eta^{-1}(\nu_a)$. 
For all $1\leq a\leq m$ the surface $S_a$ is obtained 
by succesively blowing up the total space of 
\[
M|_{C_a}\simeq K_{C_a}(n_ap_a+\nu_a)
\]
and removing a certain divisor, as shown in detail in Section \ref{toract}. The divisor in question is chosen such that the pair $(S_a, F_a)$ is a log Calabi-Yau 
surface, $K_{S_a} \simeq \CO_{S_a}(-F_a)$, for all $1\leq a\leq m$. 

Finally, let $\oY_\uxi$ be the total space of the dualizing line bundle of the normal crossing surface 
$\oS_\uxi$. 
Note that there is a natural projection 
$\pi:\oY_\uxi\to \oC$ and $\oY_\uxi$ has normal crossing singularities along the divisors $\Delta_a = \pi^{-1}(\nu_a)$, $1\leq a\leq m$.  The irreducible components of $\oY_\uxi$ are 
$Y_{a}=\pi^{-1}(C_a)$, $1\leq  a \leq m$, and $Y_0=\pi^{-1}(C_0)$.
Each $Y_a$ 
is isomorphic to the total space of the line bundle $K_{S_a}(F_a)\simeq \CO_{S_a}$ 
over $S_a$, hence $K_{Y_a}\simeq \CO_{Y_a}(-\Delta_a)$ for all $1\leq a \leq m$.
The central component $Y_0$ is isomorphic to the total space of 
the rank two bundle 
\be\label{eq:rktwocentral}
\CO_{C_0}\oplus M|_{C_0} \simeq \CO_{C_0}\oplus K_{C_0}(\nu_1+\cdots + \nu_{m}).
\ee
Hence $K_{Y_0}\simeq \CO_{Y_a}(-\Delta_1-\cdots -\Delta_{m})$ as expected.  
The log Calabi-Yau threefolds $(Y_a, \Delta_a)$ will be called wild caps in the following.

\subsection{Torus action and invariant curves}\label{toract}
A torus action on the degenerate threefold $\oY_\uxi$ will be constructed in this section satisfying the following criteria: 
\begin{itemize} 
\item  The restriction to $Y_a$, $1\leq a\leq m$, is a lift of the natural action on the rational curve $C_a$ which preserves the 
log Calabi-Yau structure. 
\item The restriction to $Y_0$ is an anti-diagonal torus action 
on the total space of the rank two bundle \eqref{eq:rktwocentral}.
\end{itemize}
The torus action will be obtained by gluing separate torus actions on each irreducible 
component of $\oY_\uxi$. 
To fix notation, let $M_a$ denote the total space of $M|_{C_a}$ and let
$(z_a, w_a)$ be standard affine coordinates on the rational component $C_a$ centered at $p_a, \nu_a$ respectively. 
Let $k_a= {\rm deg}(M|_{C_a})=n_a-1$, 
$1\leq a\leq m$. 
The inverse image of the standard open cover of $C_a$ 
determines an affine open cover $(U_a, V_a)$ of $M_a$
with coordinates $(z_a, u_a), (w_a, v_a)$ related by the transition function
\[ 
w_a=z_a^{-1}, \qquad v_a = z_a^{k_a} u_a
\]
on the overlap. Then there is a one dimensional torus action ${\bf T} \times M_a\to M_a$, 
\be\label{eq:toractMa}
t\times (z_a,u_a) \mapsto (t^{-1}z_a, u_a) \quad {\rm and}\quad 
t\times (w_a,v_a) \mapsto (tw_a, t^{k_a}v_a), 
\ee
on $M_a$ leaving the fiber over $p_a$ pointwise fixed.
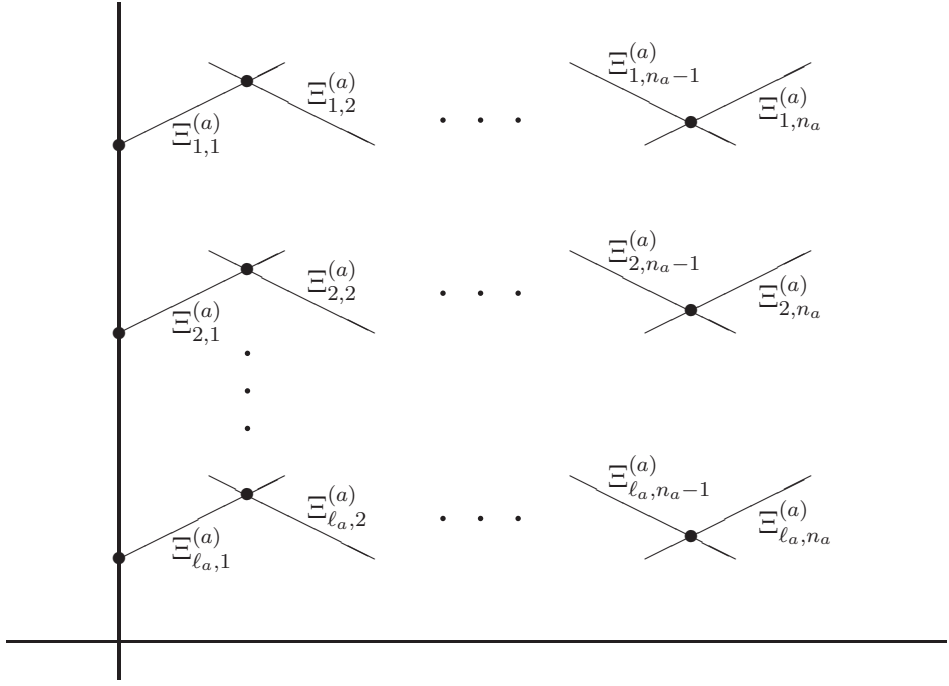
\begin{figure}
\setlength{\unitlength}{1mm}
\hspace{-20pt}\vspace{-30pt}
\setlength{\unitlength}{1mm}
\begin{picture}(170,105)
{\thicklines\put(15,10){\line(1,0){125}}
\put(30,5){\line(0,1){90}}}
\put(29,20){$\bullet$}
\put(29,50){$\bullet$}
\put(29,75){$\bullet$}
\put(30,21){\line(2,1){22}}
\put(46,28.5){$\bullet$}
\put(42,32){\line(2,-1){22}}
\multiput(72,26)(5,0){3}{\huge .}
\put(90,32){\line(2,-1){22}}
\put(105,23){$\bullet$}
\put(100,21){\line(2,1){22}}
\put(30,51){\line(2,1){22}}
\put(46,58.5){$\bullet$}
\put(42,62){\line(2,-1){22}}
\multiput(72,56)(5,0){3}{\huge .}
\put(90,62){\line(2,-1){22}}
\put(105,53){$\bullet$}
\put(100,51){\line(2,1){22}}
\put(30,76){\line(2,1){22}}
\put(46,83.5){$\bullet$}
\put(42,87){\line(2,-1){22}}
\multiput(72,79)(5,0){3}{\huge .}
\put(90,87){\line(2,-1){22}}
\put(105,78){$\bullet$}
\put(100,76){\line(2,1){22}}
\multiput(46,38)(0,5){3}{\huge .}
\put(37,21){$\Xi_{\ell_a,1}^{(a)}$}
\put(37,51){$\Xi_{2,1}^{(a)}$}
\put(37,76){$\Xi_{1,1}^{(a)}$}
\put(55,26.5){$\Xi_{\ell_a,2}^{(a)}$}
\put(55,56.5){$\Xi_{2,2}^{(a)}$}
\put(55,81.5){$\Xi_{1,2}^{(a)}$}
\put(95,30.5){$\Xi_{\ell_a,n_a-1}^{(a)}$}
\put(95,60.5){$\Xi_{2,n_a-1}^{(a)}$}
\put(95,85.5){$\Xi_{1,n_a-1}^{(a)}$}
\put(115,24.5){$\Xi_{\ell_a,n_a}^{(a)}$}
\put(115,54.5){$\Xi_{2,n_a}^{(a)}$}
\put(115,79.5){$\Xi_{1,n_a}^{(a)}$}
\end{picture}
\caption{Chains of exceptional curves on the blow-up.}
\label{SurfaceI}
\end{figure}

In order for this action to lift to the blow-up, the local sections $\xi_{a,i}: D_a \to M_{D_a}$ must be chosen to be equivariant. This means the image of each section $\xi_{a,i}:D_a\to M_a$ will be a
non-reduced zero dimensional subscheme $\delta_{a,i}\subset M_a$ given by
\[ 
z_a^{n_a}=0, \qquad u_a = \lambda_{a,i},
\]
with $\lambda_{a,i}\in \IC$, $1\leq i \leq \ell_a$. By the genericity assumption, $(\lambda_{a,i})$, $1\leq i \leq \ell_a$ must be pairwise distinct and nonzero for each $1\leq a\leq m$. 
Then there is a canonical lift of the action \eqref{eq:toractMa} to the blow-up 
along the union $\cup_{i=1}^{\ell_a} \delta_i$. 

In more detail, the surface obtained by blowing up $M_a$ 
contains $\ell_a$ linear chains of exceptional divisors 
$\Xi_{i,j}^{(a)}$, $1\leq i \leq \ell_a$, $1\leq j \leq n_a$ as shown 
in Figure 1. The surface $S_a$ is the complement of the union 
\[
\cup_{i=1}^{\ell_a} \cup_{j=1}^{n_a-1} \Xi^{(a)}_{i,j}. 
\]
Note in particular that the last component $\Xi^{(a)}_{i, n_a}$, $1\leq i\leq \ell_a$  of each 
chain is not removed, but it restricts to a non-compact curve 
on $S_a$. Moreover, the torus action \eqref{eq:toractMa} restricts to 
a torus action on $S_a$. 

By construction there is an affine open subset of $S_a$ which is canonically identified with affine chart $V_a\subset M_a$ by the blow-up map. Therefore one can use $(w_a,v_a)$ as local coordinates on $S_a$ as 
well. The threefold $Y_a$ is the total space of the line bundle $K_{S_a}(F_a)\simeq \CO_{S_a}$, where $F_a$ is the divisor $w_a=0$. In particular, $Y_a$ is a log Calabi-Yau threefold with $K_{Y_a} \simeq \CO_{Y_a}(-\Delta_a)$, where $\Delta_a$ is 
the inverse image of $F_a$. The torus action on $S_a$ lifts to a torus 
action on $Y_a$ preserving the log Calabi-Yau structure. 
In order to write a local formula, note that the inverse image of $V_a$ in $Y_a$ is an affine coordinate chart on $Y_a$ with coordinates $(w_a,y_a,v_a)$ where $y_a$ is a natural vertical affine coordinate. In this chart the torus action reads 
\be\label{eq:toractYa}
t\times (w_a, y_a,v_a) = 
(t w_a,  t^{-k_a} y_a,t^{k_a} v_a). 
\ee

It is important to note that there is a finite collection of smooth rational torus invariant curves 
$X_{a,i}$ in $Y_a$ locally given by 
\be\label{eq:invcurves} 
y_a=0, \qquad v_a = \lambda_{a, i} w_a^{k_a},\qquad 1\leq i \leq \ell_a. 
\ee
This is sketched in Figure 2. 
Each of these curves is a $(0,-1)$ curve on $Y_a$ intersecting the 
exceptional divisor $\Xi^{(a)}_{i, n_a}$ transversely at a torus fixed 
point $p_{a,i}$. In fact for each $1\leq i \leq \ell$ the blow-up construction yields an affine open coordinate chart on 
$Y_a$ with coordinates $(z_{i,a}, x_{i,a},u_{i,a})$ 
 centered at $p_{a,i}$ which are related to $(w_a,y_a,v_a)$ by the 
transition functions
\[ 
w_a=z_{i,a}^{-1}, \qquad y_a = x_{i,a},\qquad v_a = z_{i,a} u_{i,a}. 
\]
In this coordinate chart the curve $X_{i,a}$ is cut by $u_{i,a}=x_{i,a}=0$. 
Furthermore, the local form of the torus action reads 
\be\label{eq:toractYaa}
t\times (z_{i,a}, x_{i,a},u_{i,a})\mapsto (t^{-1}z_{i,a}, t^{-k_a} x_{i,a}, t^{k_a+1}u_{i,a}).
\ee

\begin{figure}
\setlength{\unitlength}{1mm}
\hspace{-70pt}\vspace{-30pt}
\begin{picture}(250,100)
{\thicklines\put(15,10){\line(1,0){180}}
\put(30,5){\line(0,1){90}}}
\put(153,7){${\bf o}$}
\put(153,86){$F_a$}
\put(152,5){\color{red}\line(0,1){90}}
\put(29,20){$\bullet$}
\put(29,50){$\bullet$}
\put(29,75){$\bullet$}
\put(30,21){\line(2,1){22}}
\put(46,28.5){$\bullet$}
\put(42,32){\line(2,-1){22}}
\multiput(72,26)(5,0){3}{\huge .}
\put(90,32){\line(2,-1){22}}
\put(105,23){$\bullet$}
\put(100,21){\line(2,1){22}}\put(122,31.5){$\bullet$}
{\color{blue}\qbezier(123,32)(152,-12)(183,32)}
\put(184,32){$X_{a,\ell_a}$}
\put(30,51){\line(2,1){22}}
\put(46,58.5){$\bullet$}
\put(42,62){\line(2,-1){22}}
\multiput(72,56)(5,0){3}{\huge .}
\put(90,62){\line(2,-1){22}}
\put(105,53){$\bullet$}
\put(100,51){\line(2,1){22}}\put(122,61.5){$\bullet$}
{\color{blue}\qbezier(123,62)(152,-42)(183,62)}
\put(184,62){$X_{a,2}$}
\put(30,76){\line(2,1){22}}
\put(46,83.5){$\bullet$}
\put(42,87){\line(2,-1){22}}
\multiput(72,79)(5,0){3}{\huge .}
\put(90,87){\line(2,-1){22}}
\put(105,78){$\bullet$}
\put(100,76){\line(2,1){22}}\put(122,86.5){$\bullet$}
{\color{blue}\qbezier(123,88)(152,-68)(183,88)}
\put(184,88){$X_{a,1}$}
\multiput(46,38)(0,5){3}{\huge .}
\put(37,21){$\Xi_{\ell_a,1}$}
\put(37,51){$\Xi_{2,1}$}
\put(37,76){$\Xi_{1,1}$}
\put(55,26.5){$\Xi_{\ell_a,2}$}
\put(55,56.5){$\Xi_{2,2}$}
\put(55,81.5){$\Xi_{1,2}$}
\put(95,30.5){$\Xi_{\ell_a,n-1}$}
\put(95,60.5){$\Xi_{2,n-1}$}
\put(95,85.5){$\Xi_{\ell_a, n-1}$}
\put(115,24.5){$\Xi_{\ell_a,n}$}
\put(115,54.5){$\Xi_{2,n}$}
\put(115,79.5){$\Xi_{1,n}$}
\multiput(184.5,42)(0,5){3}{\huge .}
\end{picture}
\caption{Torus invariant curves in $S_a$.}
\label{SurfaceII}
\end{figure}
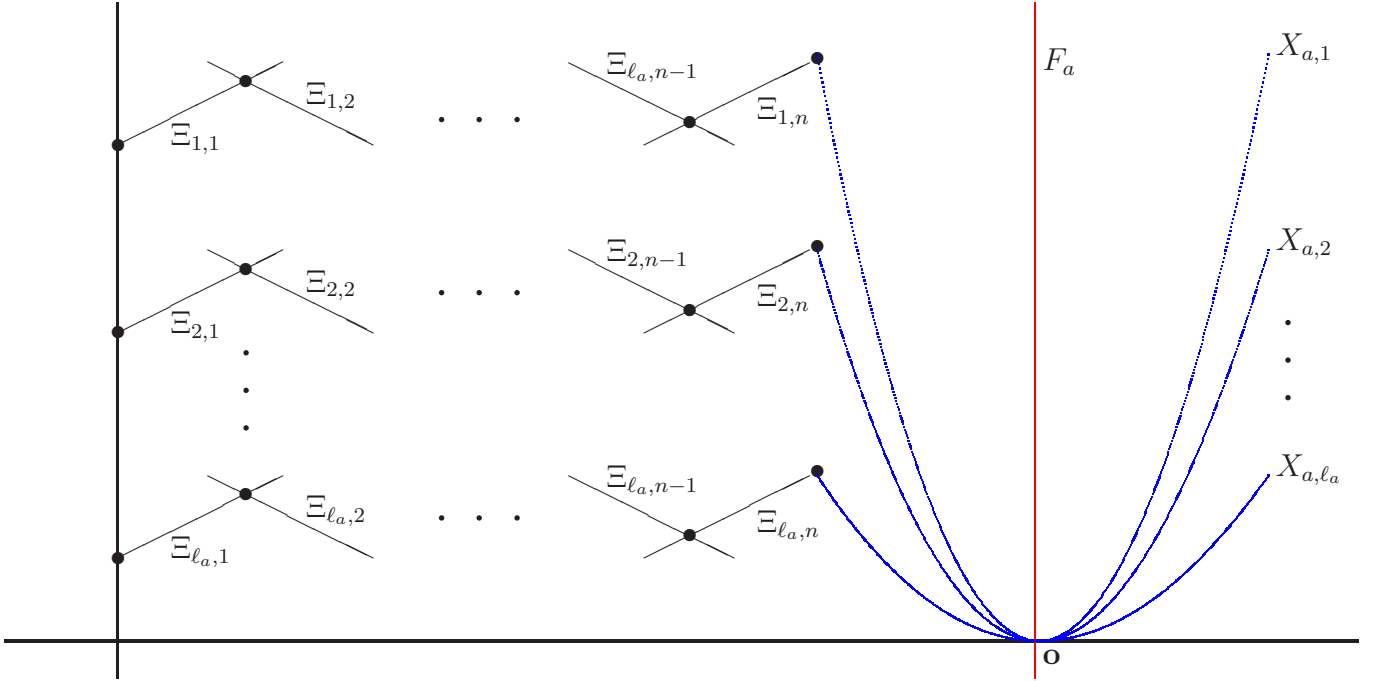

In conclusion, if the local sections $\uxi_a$ are equivariant, each irreducible component $Y_a$, $1\leq a\leq m$ has a torus action 
preserving the log Calabi-Yau structure. According to  equation \eqref{eq:toractYa}, the restriction of this action to the divisor $\Delta_a\simeq \IA^2$ has weights 
$(k_a,-k_a)$ on the tangent space at the origin, which 
is spanned by $({\partial/\partial y_a},\  {\partial/ \partial v_a})$. In order to construct a torus action on the normal crossing threefold $\oY_\uxi$, 
the weights $(k_a, -k_a)$ must be identical for all values of $a$, 
that is $k_a= k$ for all $1\leq a \leq m$.  If this condition is satisfied, there is a fiberwise torus with weights $(k,-k)$ on the central component $Y_0$ which agrees with the torus actions on 
the components $Y_a$ along the gluing divisor $\Delta_a$ for 
$1\leq a\leq m$.  This will be assumed to be the case in the next section. The log Calabi-Yau threefolds 
$(Y_{a}, \Delta_a)$, $1\leq a\leq m$ equipped with the above 
torus action will be called  degree $k$ wild caps 
in the following.

\section{Wild degenerate Gromov-Witten theory}\label{GWwildlocal} 

The goal of this section is to derive an explicit formula for the Gromov-Witten theory of a degenerate wild local curve. 
This will be carried out using J. Li's theory of relative stable morphisms 
\cite{Stable_sing} and degeneration formula \cite{Deg_formula}, the relative virtual localization theorem of Graber and Vakil \cite{Relative_loc},  as well as the local TQFT formalism of Bryan and Pandharipande \cite{Loc_curves}. Several 
results proven in \cite{MV_formula,Hodge_unknot,Curves_TQFT,GW_completed,Virasoro_curves,
Two_partition} will be also needed in the process. 

The main building blocks will be residual relative Gromov-Witten invariants of log Calabi-Yau threefolds $(Y,\Delta)$ equipped with a torus action, where $\Delta \subset Y$ is a smooth divisor. 
For completeness a review 
of relative stable maps is provided in Section \ref{relativemaps}. 
As explained there, simple relative stable maps to the pair $(Y,\Delta)$ are 
stable maps $f:\Sigma \to Y$ with prescribed contact conditions 
along $\Delta$ specified by a Young diagram 
$\lambda=(\lambda_1, \ldots, \lambda_{l_\lambda})$. More precisely, one requires that $f^{-1}(\Delta) = \sum_{j=1}^{l_\lambda} \lambda_j\sigma_j$ 
for some pairwise distinct smooth points $\sigma_1, \ldots, \sigma_{l_\lambda} \in \Sigma$. The topological invariants of a relative stable map are the 
genus $h\in \IZ$ and the homology class $\beta = f_*[\Sigma]\in H_2(Y,\IZ)$. 
As in \cite{Curves_TQFT,GW_completed,Virasoro_curves,Loc_curves}, disconnected domains are allowed 
as long as no connected component is mapped to a point.
Moreover, following loc. cit.,
the points $\sigma_1, \ldots, \sigma_{l_\lambda}$ in the domain will be unmarked, as opposed to \cite{Stable_sing,Deg_formula}.

\subsection{The wild cap}\label{wildcap}
A wild cap is a log Calabi-Yau threefold $(Y_a, \Delta_a)$, as constructed in the previous section, equipped with a torus action as in Section \ref{toract}. 
Since the relative target will be fixed throughout this subsection,
it will be simply denoted by $(Y, \Delta)$.

The first observation is that the cone of compact curve classes on $Y$ 
is freely generated by the classes $[X_i]$, $1\leq i \leq \ell$, of torus invariant curves determined by the equations \eqref{eq:invcurves}. Again the subscript $a$ will be 
dropped. 
Given any collection of non-negative integers $\um=(m_1, \ldots, m_{\ell})$, 
let $\beta(\um) = \sum_{i=1}^\ell m_i [X_{i}]$ be the corresponding curve
class, and let $\om^\bullet_{h,\um}(\lambda)$ be the 
moduli stack of genus $h$ relative stable maps to $(Y,\Delta)$
with homology class $\beta(\um)$ and relative conditions  specified by $\lambda$. 
The torus fixed locus $ \om^\bullet_{h,\um}(\lambda)^{\bf T}$ is compact and the residual relative invariants are defined 
by virtual integration 
\be\label{eq:wildcapA}
GW^\bullet({h,\um},\lambda)  = \int_{[\om^\bullet_{h,\um}(\lambda)^{\bf T}]^{vir}} {1\over e_{\bf T}(N^{vir})}
\ee
where $N^{vir}$ denotes the 
virtual normal bundle to the fixed locus. 
The wild cap is 
the generating function 
\be\label{eq:wildcapB} 
W_\lambda({\sf x}, g_s) = \sum_{h\in \IZ} \sum_{(m_1, \ldots, m_\ell) 
} g_s^{2h-2+l(\mu)} GW^\bullet(h,\um,\lambda)  \prod_{i=1}^\ell x^{m_i} 
\ee

The computation of the residual invariants \eqref{eq:wildcapA} 
proceeds by analogy with the $(0,-1)$ cap 
reviewed in Section \ref{logCYcap}. 
Using the relative virtual localization theorem \cite{Relative_loc}, the partition function of the 
 wild cap is a convolution 
\be\label{eq:wildcapC} 
W_\lambda({\sf x}, g_s) = S_\lambda({\sf x}, g_s)+\sum_{|\rho|=|\lambda|} (ik\sfh)^{2l(\rho)}\zeta(\rho) S_{\rho}({\sf x}, g_s) R_{\rho, \lambda}(g_s)
\ee
where 
\[ 
S_\rho ({\sf x}, g_s)  = \sum_{h \in \IZ} \sum_{(m_1, \ldots, m_\ell) 
} g_s^{2h-2+l(\rho)} GW^\bullet_s({h,\um},\rho)  \prod_{i=1}^\ell x^{m_i} 
\]
is the generating function of equivariant simple residues and $R_{\rho, \lambda}(g_s)$ is the generating function of rubber integrals. In the above formula  
$\sfh \in H^2(B{\bf T})$ is the equivariant parameter
and the combinatorial factor $\zeta(\rho)$ is defined by 
\[
\zeta(\rho) = \prod_{j\geq 1} k_j! j^{k_j} 
\]
for a partition $\rho = (1^{k_1}2^{k_2}\cdots)$. 

An important observation is that 
the rubber target occuring in the expanded degenerations 
of the wild cap is $\IP^1\times \IA^2$ with antidiagonal torus action 
of type $(k,-k)$ on $\IA^2$. 
Therefore using the results of \cite{GW_completed,Virasoro_curves,Two_partition} 
collected in Section \ref{rubberHurwitz}, 
the rubber generating function is
\be\label{eq:rubberfct}
R_{\rho,\lambda}(g_s)  = (ik\sfh)^{-(l(\rho)+l(\lambda))} \sum_{|\nu|=d}
\big(e^{-ikc(\nu)}-1\big)
{\chi^\nu(\rho)\over \zeta(\rho)} {\chi^\nu(\lambda)\over \zeta(\lambda)}.
\ee
for any Young diagrams $(\rho, \lambda)$ with $|\rho|=|\lambda|=d$.
Note that $\chi^{\nu}(\lambda)$ denotes the character of the irreducible 
symmetric group representation corresponding to $\nu$ evaluated on the conjugacy class 
corresponding to $\lambda$. 

The generating function $S_{\rho}({\sf x}, g_s)$ for the residues of the simple fixed loci is computed by localization arguments generalizing those used in 
Section \ref{logCYcap} for the $(0,-1)$ Calabi-Yau cap. 
The domain of a simple fixed map in the current setup is union
\[ 
\Sigma = \big(\cup_{i=1}^\ell \Sigma_i \big)\cup \big(\cup_{i=1}^\ell \cup_{j=1}^{k_i} \Lambda_{i,j}\big)
\]
where 
\begin{itemize} 
\item each $\Sigma_i$, $1\leq i \leq \ell$ is a possibly disconnected genus $h_i$ component mapped to the fixed point $p_i$ in $Y$, 
\item  for each $1\leq i \leq \ell$, the components $\Lambda_{i,1}, \ldots, \Lambda_{i,k_i}$ are projective lines attached to $\Sigma_i$ 
mapped in a torus invariant fashion to the curve $X_i$ in the target 
with some degrees $d_{i,1}, \ldots, d_{i, k_i}\geq 1$.
\end{itemize} 
Note that 
\[ 
h = \sum_{i=1}^\ell h_i, \qquad m_i = \sum_{j=1}^{k_i} d_{i,j}
\]
for all $1\leq i \leq \ell$. 
Moreover, 
the infinitesimal neighborhood of each target curve $X_i\subset Y$ is isomorphic to the infinitesimal neighborhood 
of the zero section in the $(0,-1)$ cap geometry discussed in detail in Section
\ref{logCYcap}.  
In conclusion a simple fixed locus is isomorphic to a product 
\[
\times_{i=1}^\ell \Gamma_{s}(h_i, \rho_i)
\]
where $\rho_i$ be the partition of $m_i$ determined by $(d_{i,j})$, $1\leq j \leq k_i$, and 
each factor $\Gamma_s(h_i,\rho_i)$ is a fixed locus 
of type $(h_i, \rho_i)$ in the $(0,-1)$ cap geometry.  
Furthermore, the Young diagram $\rho$ encoding the contact conditions 
along $\Delta$ is given by 
\[
\rho = \cup_{i=1}^\ell \rho_i.
\]
This yields 
\be\label{eq:wildcapD} 
S_\lambda({\sf x}, g_s) = \sum_{\rho_1\cup \cdots \cup \rho_\ell= \lambda} \prod_{i=1}^\ell\left(x_i^{|\rho_i|} S_{\rho_i}(g_s)\right)
\ee
where $S_\rho(g_s)$ is the contribution of simple fixed loci with relative 
conditions $\rho$ for the $(0,-1)$ cap geometry.  
As explained in Section \ref{simple}, the latter is given by the 
Marino-Vafa 
formula \cite{Framed_knots} proven in \cite{MV_formula,Hodge_unknot}. This reads 
\[
S_\rho(Y,\dl; g_s) = (ik{\sf u})^{-l(\rho)} \sum_{|\nu|=|\rho|} 
{\chi^{\nu}(\rho) \over \zeta(\rho)} V_\nu^{(k+1)}(q),
\]
where $q=e^{ig_s}$ and 
\[ 
V_\nu^{(k+1)}(q)=q^{(k+1)c(\nu)} s_\nu(\uq)=q^{kc(\nu)} s_{\nu^t}(\uq), \qquad \uq=(q^{1/2}, q^{3/2}, \ldots). 
\]
is the one leg topological vertex \cite{topvert}. 

The computation of the wild cap then proceeds by substituting 
equations \eqref{eq:wildcapD} and \eqref{eq:rubberfct} in \eqref{eq:wildcapC}. As shown in detail below, the final expression
is 
\be\label{eq:wildcapE} 
\bal 
W_\lambda({\sf x}, g_s) = 
(ik\sfh)^{-l(\lambda)} \zeta(\lambda)^{-1}
\sum_{|\mu|=|\lambda|} \sum_{\substack{ \nu_1, \ldots, \nu_\ell\\ |\nu_1|+\cdots + |\nu_\ell| = |\lambda|}} q^{-kc(\mu)}
c^\mu_{\nu_1, \ldots, \nu_\ell} \chi^\mu(\lambda) 
\prod_{i=1}^\ell\left(x_i^{|\nu_i|} q^{kc(\nu_i)} s_{\nu_i^t}(\uq)\right)
\\
\eal 
\ee
where $c^\mu_{\nu_1,\ldots, \nu_\ell}$ are the fusion coefficients for Schur functions i.e. 
\be\label{eq:fusionA} 
\prod_{i=1}^\ell s_{\nu_i}({\sf x}) = \sum_{|\mu|=
 |\nu_1|+\cdots + |\nu_\ell|}  
c^\mu_{\nu_1, \ldots, \nu_\ell} s_\mu({\sf x}).
\ee
The proof of formula \eqref{eq:wildcapE} will proceed in a few steps. 

{\bf Step 1.} Using the orthogonality relation 
\[
\sum_{|\mu|=|\lambda|=|\rho|} \chi^\mu(\rho)\chi^\mu(\lambda) = \zeta(\lambda)\delta_{\rho, \lambda},
\]
the 
 generating function $S_\lambda({\sf x}, g_s)$ in 
\eqref{eq:wildcapD} will be rewritten as follows 
\[
\bal 
S_\lambda({\sf x}, g_s) & =  (iku)^{-l(\lambda) } 
\sum_{
\rho_1\cup \cdots \cup \rho_\ell = \lambda }
\sum_{\substack{\nu_1, \ldots, \nu_\ell\\ |\nu_i|=|\rho_i|,\ 
1\leq i \leq \ell}}
\prod_{i=1}^\ell \left({x_i}^{|\rho_i|}\zeta(\rho_i)^{-1}\chi^{\nu_i}({\rho_i})V^{(k+1)}_{\nu_i}(q)\right) \\
&  =  (iku)^{-l(\lambda) } \zeta(\lambda)^{-1} \sum_{|\rho|=|\lambda|}
\sum_{\rho_1\cup \cdots \cup \rho_\ell = \rho} 
\zeta(\lambda) \delta_{\rho, \lambda} \\
&\qquad\ \qquad\qquad\
\sum_{\substack{\nu_1, \ldots, \nu_\ell\\ |\nu_i|=|\rho_i|,\ 
1\leq i \leq \ell}}
\prod_{i=1}^\ell \left({x_i}^{|\rho_i|}\zeta(\rho_i)^{-1}\chi^{\nu_i}({\rho_i})V^{(k+1)}_{\nu_i}(q)\right)\\
\eal 
\]
\[
\bal
{}
& =  (iku)^{-l(\lambda) } \zeta(\lambda)^{-1} \sum_{|\mu|=|\lambda|} \sum_{|\rho|=|\lambda|}
\sum_{\rho_1\cup \cdots \cup \rho_\ell = \rho} 
\chi^\mu(\rho)\chi^\mu(\lambda) \\
&\qquad\ \qquad\qquad\
\sum_{\substack{\nu_1, \ldots, \nu_\ell\\ |\nu_i|=|\rho_i|,\ 1\leq i \leq \ell}}
\prod_{i=1}^\ell \left({x_i}^{|\rho_i|}\zeta(\rho_i)^{-1}\chi^{\nu_i}({\rho_i})V^{(k+1)}_{\nu_i}(q)\right)\\
& =  (iku)^{-l(\lambda) } \zeta(\lambda)^{-1}
\sum_{\substack{\nu_1, \ldots, \nu_\ell\\ 
|\nu_1|+\cdots + |\nu_\ell|=|\lambda|}}
\sum_{|\mu|=|\lambda|}
\chi^\mu(\lambda) \prod_{i=1}^\ell \left(x_i^{|\nu_i|} V_{\nu_i}^{(k+1)}(q)\right)\\
& \qquad\qquad \qquad\ \sum_{\substack{\rho_1, \ldots,\rho_\ell\\ 
|\rho_i|=|\nu_i|,\ 1\leq i \leq \ell}} \chi^\mu(\rho_1\cup \cdots \cup \rho_\ell)\prod_{i=1}^\ell
\left(\zeta(\rho_i)^{-1} \chi^{\nu_i}({\rho_i})\right) 
\eal 
\]
Now note the following combinatorial formula
\be\label{eq:fusioncoeff}
c^\mu_{\nu_1, \ldots, \nu_\ell}=\sum_{\substack{\rho_1, \ldots,\rho_\ell\\ 
|\rho_i|=|\nu_i|,\ 1\leq i \leq \ell}} \chi^\mu(\rho_1\cup \cdots \cup \rho_\ell)\prod_{i=1}^\ell
\left(\zeta(\rho_i)^{-1} \chi^{\nu_i}({\rho_i})\right) 
\ee
for the fusion coefficients, which 
can be easily proven by writing the Schur functions in terms of power sum symmetric functions. This yields by substitution:
\be\label{eq:wildcapF}
\bal 
S_\lambda({\sf x}, g_s) & =  (iku)^{-l(\lambda) } 
\zeta(\lambda)^{-1}
\sum_{\substack{\nu_1, \ldots, \nu_\ell\\ |\nu_1|+\cdots + |\nu_\ell|=|\lambda|}}
\sum_{|\mu|=|\lambda|} c^\mu_{\nu_1, \ldots, \nu_\ell}\chi^\mu(\lambda) \prod_{i=1}^\ell \left(x_i^{|\nu_i|} V_{\nu_i}^{(k+1)}(q)\right). 
\eal 
\ee
for the generating function of simple residues. 

{\bf Step 2}. The second term in the convolution formula \eqref{eq:wildcapC} is computed using equation \eqref{eq:wildcapD} and identity \eqref{eq:fusioncoeff}
\be\label{eq:wildcapG}
\bal 
 \sum_{|\rho|=|\lambda|} (ik\sfh)^{2l(\rho)}\zeta(\rho) 
S_{\rho}({\sf x}, g_s) R_{\rho, \lambda}(g_s) & =  
(ik\sfh)^{-l_\lambda} \zeta(\lambda)^{-1} \sum_{|\rho|=|\lambda|} \sum_{|\mu|=|\lambda|} \big(e^{-ikc(\mu)}-1\big) \chi^\mu(\rho)\chi^\mu(\lambda)\\
& \ \ \sum_{\rho_1\cup\cdots\cup\rho_\ell=\rho} \sum_{\substack{\nu_1, \ldots, \nu_\ell\\ |\nu_i|=|\rho_i|,\ 1\leq i \leq \ell}}
\prod_{i=1}^\ell\left(x_i^{|\rho_i|}
\zeta(\rho_i)^{-1}\chi^{\nu_i}(\rho_i) V_{\nu_i}^{(k+1)}(q)\right)\\
& =  (ik\sfh)^{-l_\lambda} \zeta(\lambda)^{-1} 
\sum_{|\mu|=|\lambda|} \big(e^{-ikc(\mu)}-1\big)
 \chi^\mu(\lambda)\\
& \ \  \sum_{|\nu_1|+\cdots + |\nu_\ell|=|\lambda|} 
c^\mu_{\nu_1, \ldots, \nu_\ell} \prod_{i=1}^\ell
\left(x_i^{|\nu_i|} V_{\nu_i}^{(k+1)}(q)\right)\\
\eal 
\ee
Finally, substituting \eqref{eq:wildcapF} and \eqref{eq:wildcapG}
in \eqref{eq:wildcapC} one obtains formula \eqref{eq:wildcapE}. 

\subsection{The central component}\label{central}

The central component of the degenerate threefold $\oY_\uxi$ is isomorphic to the total space of the rank two bundle 
\[ 
\CO_{C_0}\oplus 
K_{C_0}(\nu_1+\cdots +\nu_{m})
\]
over the smooth genus $g$ curve $C_0$. Under the current assumption there is a fiberwise torus action with weights $(k,-k)$ on $Y_0$. Recall that $Y_0$ is glued to the wild caps $Y_1, \ldots, Y_m$ along thefibers $\Delta_1, \ldots, \Delta_{m}$ over the points $\nu_1, \ldots, \nu_{2m}$. 
The goal of this section is to determine the generating function of relative Gromov-Witten invariants for the 
data $(Y_0, \Delta_1, \ldots, \Delta_{m})$ with contact conditions 
$\lambda_1, \ldots, \lambda_{m}$. An explicit formula for this generating function 
follows from the two dimensional TQFT formalism for local curves constructed in  \cite{Loc_curves}. The main building blocks are determined in the proof of 
Theorems 7.3 in loc. cit. as explained below, mainly using the same 
notation. 

For the antidiagonal torus action, the TQFT constructed in \cite{Loc_curves} takes values in the category of $R$-modules, where $R=\IQ[i](\sft)((g_s))$ 
is the ring of Laurent series in the genus counting variable $g_s$ 
over the algebraic extension  $\IQ[i]$. 
Given a fixed degree $d\geq 1$, the $R$-module  ${\bf GW}_d(S^1)$
assigned to the circle is freely generated by elements $e_\alpha$ in one-to-one correspondence with Young diagrams 
$|\alpha|= d$. Let $e^\alpha$ be the dual generators given by 
\[ 
e^\alpha(e_\beta) = (-\sft^2)^{l(\alpha)} \zeta(\alpha) \delta^\alpha_\beta. 
\] 
Note that the coefficients $(-\sft^2)^{l(\alpha)} \zeta(\alpha)$ are the gluing factors in the degeneration formula for local relative Gromov-Witten invariants  
\cite[Thm. 3.2]{Loc_curves}.  
The multiplication operator is determined by the residual relative invariants 
of a local genus zero $(0,0)$ curve with contact conditions along three fibers
\be\label{eq:multop}
{\sf M}= \sum_{|\alpha|=|\beta|=|\gamma|=d} 
{\sf GW}(0|0,0)^\gamma_{\alpha,\beta} e^\alpha\otimes e^\beta \otimes e_\gamma.
\ee
The genus operator is defined as 
\be\label{eq:genusop}
{\sf G} = \sum_{|\alpha|=|\beta|=d} 
{\sf GW}(1|0,0)_\alpha^\beta e^\alpha \otimes e_\beta 
\ee
where the coefficients ${\sf GW}(1|0,0)_\alpha^\beta$ are the 
 residual relative invariants 
of a local genus one $(0,0)$ curve
one $(0,0)$ curve with 
contact conditions along two fibers. 
An important observation 
is that the multiplication operator takes the canonical form 
\be\label{eq:multopB}
{\sf M}^\mu_{\lambda, \rho} = \delta^\mu_\lambda \delta^\mu_\rho. 
\ee
with respect to the idempotent basis, 
\be\label{eq:vbasis} 
v_\lambda = \zeta(\lambda)^{-1} \sum_{|\alpha|=d} (i\sft)^{l(\alpha)-d}\chi^\lambda(\alpha) e_\alpha. 
\ee
The  remaining TQFT elements needed in this section are:
\begin{itemize}
\item the $(0,1)$ annulus in the $v$-basis,
\be\label{eq:vcyl}
\sfA^\rho_\lambda =   (i\sft)^{-d}  \zeta(\lambda)^{-1} s_{\lambda^t}(\uq)^{-1}\delta_{\rho, \lambda}, 
\ee
\item the genus operator in the $v$-basis,
\be\label{eq:genusop} 
{\sf G}_{\lambda, \rho} = (i\sft)^{2d} \zeta(\lambda)^2\delta_{\lambda, \rho},
\ee
\item and the counit in the $v$-basis, which is determined by the $(0,0)$ cap, 
\be\label{eq:counit}
{\sf C}_\lambda = (i\sft)^{-2d} \zeta(\lambda)^{-2}.
\ee
\end{itemize}

In this formalism, a genus $g$ local 
curve of type $(p,q)$ with marked points $\unu=(\nu_1, \ldots,\nu_m)$ 
determines an $R$-module morphism ${\bf GW}_d(C,\unu|p,q): {\bf GW}_d(S^1)^{\otimes m} \to R$ given by 
\[ 
\bal
{\bf GW}_d(C,\unu|p,q)& = \sum_{|\alpha_1|=\cdots = |\alpha_m|=d} {\sf GW}(g|p,q)_{\alpha_1, \ldots, \alpha_s} e^{\alpha_1}\otimes \cdots \otimes e^{\alpha_s}.
\eal
\]
Here ${\sf GW}(g|p,q)_{\alpha_1, \ldots, \alpha_s}$ are residual 
 relative Gromov-Witten invariants of the local $(p,q)$ curve $C$ with contact conditions $(\alpha_1, \ldots, \alpha_s)$ along the $m$ fibers $\Delta_1, \ldots, \Delta_m$.  
In the present case, recall that $(p,q)=(0,2g-2+m)$ and $\sft = k\sfh$. 
In order to write down an explicit formula, let 
${\sf P}^m: {\bf GW}_d(S^1)^{\otimes m} \to {\bf GW}_d(S^1)$ be the operator 
\[ 
{\sf P}^m(v_1, \ldots, v_m) = {\sf M}({\sf M}(\cdots({\sf M}(v_1,v_2),v_3),\ldots),v_m). 
\]
Then note that 
\be\label{eq:vcentralA}
{\bf GW}_d(C,\unu|0,2g-2+m) = {\sf C} {\sf A}^{2g-2+m} {\sf G} {\sf P}^m. 
\ee
Using equations \eqref{eq:multopB}, \eqref{eq:vcyl}, \eqref{eq:genusop} and \eqref{eq:counit}, the coefficients of this operator in the $v$-basis are 
\be\label{eq:vcentralB}
Z_{\lambda_1, \ldots, \lambda_m}=  (ik\sfh)^{-(2g-2+m)d} 
 \sum_{|\rho|=d} 
s_{\rho^t}(\uq)^{-(2g-2+m)} \prod_{a=1}^m \delta^\rho_{\lambda_a}. 
\ee

\subsection{The wild curve formula}\label{wildformula}

The residual Gromov-Witten of a wild degenerate local curve is obtained from equations \eqref{eq:wildcapE} and \eqref{eq:vcentralA}, \eqref{eq:vcentralB}  
using the degeneration formula. For fixed degree $r\geq 1$ each wild cap yields an element
\[
w = \sum_{|\alpha|=r} (ik\sfh)^{2l(\alpha)} \zeta(\alpha)
W_\alpha({\sf x}, g_s) e_\alpha
\]
in ${\sf GW}_d(S^1)$. The degree $r$ degenerate wild curve partition function is then given 
by 
\be\label{eq:wildcurveA}
Z_r = {\bf GW}_r(C,\unu|0,2g-2+m)(w_1\otimes \cdots \otimes w_m). 
\ee
In order to derive an explicit formula, one needs the coefficients of the wild cap in the 
$v$-basis.  Using the inverse change of basis 
\be\label{eq:invbasis}
e_\alpha = (it)^{r-l(\alpha)} \zeta(\alpha)^{-1}\sum_{|\lambda|=r} \zeta(\lambda) \chi^{\lambda}(\alpha)v_\lambda.
\ee
one obtains 
\[ 
\bal
w^{\lambda} = \zeta(\lambda)\sum_{|\alpha|=r}  (ik\sfh)^{r+l(\alpha)}
\chi^{\lambda}(\alpha) W_\alpha({\sf x}, g_s).\\
\eal 
\] 
Using equation \eqref{eq:wildcapE}, this yields 
\be\label{eq:vwildcap}
w^\lambda = (ik\sfh)^{|\lambda|} \zeta(\lambda)
q^{-kc(\lambda)} \sum_{|\nu_1|+\cdots+|\nu_\ell|=|\lambda|} 
c^\lambda_{\nu_1, \ldots, \nu_\ell} \prod_{i=1}^\ell x_i^{|\nu_i|}q^{kc(\nu_i)} s_{\nu_i^t}(\uq).
\ee
Substituting \eqref{eq:vcentralB} and \eqref{eq:vwildcap} in \eqref{eq:wildcurveA}, it follows that 
\[
\bal 
Z_r & = \sum_{|\lambda|=r} s_{\lambda^t}(\uq)^{2-2g-m} \prod_{a=1}^m F_{k,\ell_a,\lambda}({\sf x}_a,q) 
\eal
\]
where 
\[
F_{k,\ell,\lambda}({\sf x}, q) = q^{-kc(\lambda)}\sum_{|\nu_1|+\cdots+|\nu_\ell|=|\lambda|} 
c^\lambda_{\nu_1, \ldots, \nu_{\ell}} \prod_{i=1}^\ell x_i^{|\nu_i|}q^{kc(\nu_i)} s_{\nu_i^t}(\uq).
\]
This proves formula \eqref{eq:wildformulaA}.

\section{Relative stable maps}\label{relativemaps} 

This section is a self-contained review of relative stable 
maps following 
\cite{Stable_sing,Deg_formula}, and some of their applications 
\cite{MV_formula, Hodge_unknot,Curves_TQFT,Two_partition,Virasoro_curves, GW_completed, Loc_curves, Relative_loc}. No results 
presented here are new, the goal being to collect all results needed in Section \ref{wildcap} in a self-contained manner in order to streamline the derivation of the 
wild cap formula. This review is primarily aimed at the non-expert reader. 

\subsection{Background}\label{relbackground} 
 Suppose $Y$ is a smooth complex 
projective variety and $\Delta\subset Y$ is a smooth connected divisor. A simple relative stable map to the pair $(Y,\Delta)$  consists of the data $(\Sigma, f, \sigma_1, \ldots, \sigma_l )$ 
where $\Sigma$ is a connected nodal curve
 $f:\Sigma\to Y$ a map to $Y$, 
and $\sigma_1, \ldots, \sigma_l$ are smooth pairwise distinct points on $\Sigma$ such that:
\begin{itemize}
\item $f^{-1}(\dl) = \sum_{i=1}^l m_i \sigma_i$ 
for some positive integers $m_1, \ldots, m_l$, and 
\item the automorphism group of the data $(\Sigma, f, \sigma_1, \ldots, \sigma_l)$ is finite. 
\end{itemize} 
The topological invariants of simple relative stable maps are the arithmetic genus $g\geq 0$ of the domain and the homology class $\beta = f_*[\Sigma] \in H_2(X, \IZ)$. 

The moduli stack of simple relative 
stable maps with fixed invariants $(g, l, \beta)$, $(m_1, \ldots, m_l)$ is 
not proper. A compactification is constructed in \cite{Stable_sing} by allowing the target to degenerate in a controlled way. The allowed degenerations of the target are normal crossing varieties 
constructed by gluing $Y$ and an arbitrary number of copies of 
the projective bundle 
$$P = \IP_\Delta(N_{\Delta/Y}\oplus \CO_\Delta),$$
 where $N_{\Delta/Y}$ is the normal bundle of $\Delta$ in $Y$. This projective bundle has two canonical sections $\Delta_0, \Delta_\infty$ with normal bundles $N_{\dl/Y}^{-1}$, $N_{\dl/Y}$ respectively. For any integer $n\geq 1$ 
let $P_n$ be the normal crossing variety obtained by gluing $n$ copies of $P$ such that the section $\dl_\infty$ of the $i$-th copy is identified 
with section $\dl_0$ of the $(i+1)$-th copy. Note that no gluing occurs along $\dl_0$ in the first copy and $\dl_\infty$ in the last. 
Abusing notation, these 
two divisors on $P_n$ will be denoted by
$\dl_0, \dl_n$ respectively. The singular locus of $P_n$ consists of the union $\dl_1\cup \cdots \cup\dl_{n-1}$ of copies of $\dl$ where 
$\dl_i$ is the intersection between the $i$-th and the $(i+1)$-th copy of $P$. For further reference let ${\rm Aut}_{\dl}(P_n)\simeq (\IC^\times)^n$ denote the group of automorphisms of $P_n$ 
acting trivially on all sections $\Delta_0, \ldots, \dl_n$ over $\Delta$. 
The degenerate targets $Y_n$ are constructed by gluing $Y$ to $P_n$ so that $\dl\subset Y$ is identified with $\dl_0\subset P_n$. 
Therefore one obtains a relative pair $(Y_n, \dl_n)$ such that the  singular set of $Y_n$ is the union of the $n$ 
copies $\dl_0\cup  \cdots\cup \dl_{n-1}$ of copies of  $\dl$. 
Note that there is an algebraic stack ${\mathfrak E}$ of expanded 
degenerations parameterizing all degenerate targets which occur in this construction, and a universal family $\CY \to {\mathfrak E}$. 

Compactification of the moduli stack of relative stable maps is achieved by including relative stable morphisms to degenerate pairs $(Y_n, \dl_n)$. Such relative maps consist of data $(\Sigma, f, \sigma_1, \ldots, \sigma_l)$ satisfying the following conditions 
\begin{itemize} 
\item $\Sigma$ is a connected nodal curve and $f:\Sigma \to Y_n$ is a predeformable map to $Y_n$. 
\item $\sigma_1, \ldots, \sigma_l$ are smooth points of $\Sigma$ such that 
\[ 
f^{-1}(\dl_n) = \sum_{i=1}^l m_i \sigma_i.
\]
\end{itemize} 
Predeformable morphisms $f:\Sigma \to Y_n$ are defined by the following conditions:
\begin{itemize} 
\item The set theoretic inverse image of any gluing divisor $\dl_i\subset Y_n^{sing}$ is a set of nodes of $\Sigma$. 
\item The two branches of $\Sigma$ crossing at each such node are mapped to different irreducible components of $Y_n$ meeting along $\dl_i$ such the contact orders along $\dl_i$ are equal. 
\end{itemize}
The topological invariants of a relative map to an expanded target are the arithmetic 
genus $g$ of the domain, and the homology class $\beta=\pi_{n*}f_*[\Sigma]\in H_2(Y,\IZ)$ where $\pi_n : Y_n \to Y$ is the natural projection. 

Two relative maps $(\Sigma, f)$, $(\Sigma', f')$ to $Y_n$ are isomorphic if there is an isomorphism $\psi: \Sigma 
\to \Sigma'$ and an automorphism $\varphi:Y_n \to Y_n$ fixing 
$Y$ and the divisors $\dl_0, \ldots, \dl_{n}$ 
pointwise, such that $f' \circ \psi = \varphi \circ f$. Note that  $\varphi|_{P_n}\in {\rm Aut}_\dl(P_n)$ acts by fiberwise scalar multiplication 
on each copy of $P$ used in the construction of $P_n$. 
Stability is defined by requiring the relative maps to have a finite automorphism group 
according to this notion of isomorphism.

The main result of \cite{Stable_sing} proves that for fixed data $(g,l,\beta)$, 
$\um=(m_1, \ldots, m_l)$ there is a proper Deligne-Mumford moduli stack 
$\om_{g,\beta}(Y,\dl, \um)$ of relative stable morphisms equipped with a perfect obstruction theory and a universal relative morphism to the universal family $\CY$ of expanded degenerations. 

A slightly different flavor of relative theory will be 
used throughout this paper, as in \cite{Curves_TQFT,GW_completed, Virasoro_curves,Loc_curves}. Namely, as opposed to \cite{Stable_sing},  the points $(\sigma_1,\ldots, \sigma_l)$ in the inverse image 
of the divisor $D$ are unmarked, and one also allows 
disconnected domains as long as no connected component is contracted to a point. In particular the numerical invariants 
$(m_1, \ldots, m_l)$ are unordered, hence they are encoded in a Young diagram $\mu$ with $l$ rows. The moduli stack of stable relative maps with fixed invariants $(g,\beta)$,  $\mu$ will  be 
denoted by $\om_{g,\beta}(Y,\dl;\mu)$ for connected domains, 
respectively $\om_{g,\beta}^\bullet(Y,\dl;\mu)$ for disconnected 
domains. Note that in the second case $g$ is allowed to be negative. 
In both cases the moduli stack is equipped with a virtual cycle and a universal 
morphism to the stack of expanded degenerations $\CY$. 

Finally, given a torus action on $Y$ which preserves $\dl$, 
there is a virtual localization formula for moduli stacks 
of relative stable maps, which has been extensively used in the 
literature. The foundations have been proven in \cite{Relative_loc}. 
In particular, this allows one to define residual relative invariants 
for non-compact pairs $(Y,\dl)$ provided that there is a torus 
action on $Y$ which preserves $\dl$ and has compact fixed locus.  

\subsection{Rubber target}\label{rubberbackground}
The classification of torus fixed loci in moduli stacks of relative 
maps leads naturally to a variant of the above construction
employing non-rigid targets. Namely one constructs a moduli stack of 
predeformable stable maps to degenerate targets $P_n$ 
imposing fixed contact conditions along the sections $\dl_0, \dl_\infty$ specified by two partitions $\rho, \mu$. 
Moreover two such maps  
$(\Sigma, f)$, $(\Sigma', f')$ are isomorphic if there is an isomorphism $\psi: \Sigma 
\to \Sigma'$ and an automorphism $\varphi\in {\rm Aut}_\dl(P_n)$ 
such that $f' \circ \psi = \varphi \circ f$. 
This construction yields a Deligne-Mumford moduli stack $\om_{h,\beta}^\bullet(P, \rho, \mu)_\sim$ equipped with a perfect obstruction theory, where $\beta\in H_2(P,\IZ)$.  Again, disconnected domains are allowed. 

In this case there is again a stack of expanded degenerations ${\mathfrak E}$, a universal family $\calP\to 
{\mathfrak E}$, and a universal map to $\calP$. Moreover, 
${\mathfrak E}$ is naturally isomorphic to an open substack of the algebraic moduli stack  ${\mathfrak M}_{0,2}$ of 
genus zero curves with two marked points $0, \infty$. 
In particular one obtains a cohomology class $\psi_0$ on the rubber moduli stack by pulling back the $\psi$-class associated to the marked point $0$ on ${\mathfrak M}_{0,2}$. Abusing the language, 
the class $\psi_0$ will be referred to as the rubber 
$\psi$-class at $0$.

\subsection{Local $(0,-1)$ rational curves}\label{logCYcap} 

As an example, this section presents the computation of relative 
Gromov-Witten theory by localization for a local $(-1,0)$ rational 
curve. Although this computation has already been carried out for example in \cite{MV_formula, Hodge_unknot,Curves_TQFT}, a self-contained review will be helpful since many details are needed in  Section \ref{wildcap}. 

Let $C=\IP^1$ and let $(U,z)$, $(V,w)$ denote the standard affine coordinate charts on $C$. 
Let ${\bf T}\times C \to C$ be the  torus action on $C$ 
given by 
\[ 
(t,z) \mapsto t^{-1}z, \qquad (t,w)\mapsto tw
\]
in the two charts respectively. Let also $p\in U$ and $\delta \in V$ be  the points $z=0$, $w=0$ respectively. 
Let $Y$ be the total space of $\CO_C\oplus \CO_C(-p)$ and let  $\Delta$ be the fiber of $Y$ over $\delta$. 
The open subsets $Y_U, Y_V$ are affine coordinate charts on $Y$ with coordinates $(z,u,x)$, $(w,v,y)$ related by the transition functions 
\[ 
w=z^{-1}, \qquad x=y, \qquad  v = zu.
\]
Let ${\bf T}\times Y\to Y$ denote the lift of the torus action to $Y$ such that 
\[
t\times (z,u,x) \mapsto (t^{-1}z, t^{-k}x, t^{k+1}u), \qquad 
t\times (w,v,y) \mapsto (tw, t^{-k}y, t^k v).
\]

Note that in this case any expanded degeneration $Y_n$, $n\geq $,  is  naturally isomorphic to the total space of the rank two bundle $\CO_{C_n}\oplus \CO_{C_n}(-p)$ where $C_n$ is the $n$ step degeneration of the  target curve. As explained in Section \ref{relbackground}, $C_n$ is constructed by gluing $C$ to a linear chain of $n$ rational curves, $(\IP^1)_n$, where each curve in chain contains two marked points $\delta_0, \delta_\infty$. 
Similarly, any rubber target 
$P_n$ is isomorphic to $(\IP^1)_n \times \IA^2$ and the two divisors 
$\Delta_0, \Delta_n$ are the fibers of the projection to $(\IP^1)_n$
over the marked points $\delta_0, \delta_n$. 
The above torus action lifts to a torus action 
${\bf T}\times Y_n\to Y_n$ on any expanded degeneration of the target which scales the fibers over $C_n$ with weights $(k,-k)$ leaving the zero section pointwise fixed. This holds for $P_n$ as well.

The residual invariants of the pair $(Y, \Delta)$ are defined by 
evaluating an equivariant obstruction class on the virtual cycle of the moduli space of relative stable 
maps to $(C, \infty)$. Clearly, the degree of any such relative map is given  by $d= \sum_{i=1}^l \mu_i = |\mu|$, where $\mu=(\mu_1, \ldots, \mu_l)$ 
is the partition specifying the relative conditions at $\delta_\infty$.  Therefore the moduli stack of 
relative stable maps with disconnected domains can be denoted by $\om^\bullet_h(C,\delta; \mu)$, where $h\in \IZ$ is the arithmetic genus of the domain. 

Let $\phi:\CS \to \CC$ denote the universal relative stable map to $(C,\delta_\infty)$ and $\pi: \CS\to \om^\bullet_h(C,\delta;\mu)$ denote the natural projection. Note that there is a line bundle $\CO_\CC(-\wp)$ which restricts to $\CO_{C_n}(-p)$ on each 
closed point. Then the obstruction is the equivariant $K$-theory class 
\[ 
Ob = -T^{-k} R\pi_*\phi^* \CO_\CC(-\wp) - T^{k} R\pi_*\phi^*\CO_C
\]
where $T$ denotes the canonical representation of the torus and 
$\CO_{\CC}(-\wp)$, $\CO_\CC$ are equipped with their natural equivariant structures.  The relative local invariants are defined by 
\be\label{eq:localrelA}
GW^\bullet(h, \mu) = \int_{[\om^\bullet_h(C,\delta; \mu)]^{vir}_{\bf T}} e_{\bf T} (Ob), 
\ee
where $[\om^\bullet_h(C,\delta; \mu)]^{vir}_{\bf T}$ is the equivariant virtual cycle. 
They can be computed by relative virtual localization as explained below.

The torus fixed loci in  moduli stacks of relative stable maps 
have been analyzed in detail for example in \cite{MV_formula,Hodge_unknot, Curves_TQFT,Two_partition}.  
The main observation is that the domain of a generic torus invariant relative map $f:\Sigma \to C_n$ splits naturally as a union of two curves
$\Sigma = \Sigma_s \cup \Sigma_\infty$ whose intersection is 
a finite set of nodes of $\Sigma$. The 
restriction $f_s=f|_{\Sigma_s}$
maps $\Sigma_s$ to $C$ while $f_\infty = f|_{\Sigma_\infty}$ maps 
$\Sigma_\infty$ to $(\IP^1)_n$. Using the predeformability condition, it follows that $f_s:\Sigma_s\to C$ must be a simple relative stable map to $(C, \delta)$ with relative invariants specified by a Young diagram $\rho$ such that $|\rho|=|\mu|$. At the same time 
$f_\infty : \Sigma_\infty \to (\IP^1)_n$ is a relative stable map to the triple $((\IP^1)_n, \delta_0, \delta_n)$ with relative invariants $(\rho, \mu)$ respectively. Using the definition of isomorphisms of relative maps, 
it follows that the data $(\Sigma_\infty, f|_{\Sigma_\infty})$ determines a point in the rubber moduli stack of stable relative maps to $\IP^1$. 

At the same time, as explained in 
\cite[Appendix A]{Curves_TQFT}, identifying the ramification divisors on the two components $\Sigma_s$, 
$\Sigma_\infty$ over $\delta=\delta_0$ requires ordering
the ramification points in the two components, which are 
unmarked in the present construction. This implies that  a generic fixed locus $\Gamma(h,\mu)$ in the moduli stack $\om_{h}^\bullet(C,\delta, \mu)$ is isomorphic to  a finite etale cover of the product 
$\Gamma_s(h_1, \rho)\times \om^\bullet_{h_2}(\IP^1,\delta_0,\delta_\infty; \rho, \mu)$ 
where $h=h_1+h_2$ while $\rho$ is a Young diagram with 
$|\rho|=|\mu|$ and $\Gamma_s(h_1, \rho)$ is a simple fixed locus in the moduli stack 
$\om_{h_1}^\bullet(C,\delta, \rho)$.

Using the relative virtual localization theorem \cite{Relative_loc}, the equivariant residue of such a fixed locus 
factors as 
\be\label{eq:factorA}
\bal 
Z_{\Gamma(h,\mu)} & = \zeta(\rho)e_{\bf T}(T_\delta \Delta_0) Z_{\Gamma_s(h_1,\rho)}
\int_{[\om^\bullet_{h_2}(\IP^1,\delta_0,\delta_\infty; \rho, \mu) ]^{vir}} {e_{\bf T}(Ob_\infty)\over -{\sf u}- \psi_0}\\
\eal 
\ee
where 
\begin{itemize}
\item $Z_{\Gamma_s(h_1,\rho)}$ is the equivariant residue of the simple fixed locus $\Gamma(h_1, \rho)$,
\item $Ob_\infty$ is an equivariant obstruction class on the rubber moduli space, 
\item $\psi_0$ is the rubber $\psi$-class at 0 defined as in Section 
\ref{rubberbackground},
\item ${\sf u} \in H^2_{\bf T}({\rm point})$ is the natural generator of the equivariant cohomology ring of the point, and 
\item given a partition $\rho=(1^{k_1}2^{k_2} \cdots)$, the factor $\zeta(\rho)$ is defined by
$\zeta(\rho)= \prod_{j\geq 1} k_j! j^{k_j}$. 
\end{itemize}
Note that the factor $1/(-{\sf u} - \psi_0)$ corresponds to normal 
infinitesimal deformations in the ambient moduli stack $\om_h^\bullet(C,\delta; \mu)$ induced by deformations of the degenerate target. As shown in \cite[Eqn. 14]{Curves_TQFT}, the combinatorial factor $\zeta(\rho)$ encodes the degree of the finite etale cover involved 
in the presentation of the fixed locus as a direct product.  

The rubber obstruction class is given by 
\[ 
Ob_\infty = \left(-T^k [R\pi_*\phi^*\CO_{\calP^1}] -
T^{-k} [R\pi_*\phi^*\CO_{\calP^1}]\right)\big|_{\Gamma_\infty(h,\rho,\mu)}
\]
where $\phi: \CS \to \calP^1$ is the universal relative morphism 
to the rubber target and $\pi: \CS \to \om^\bullet_{h_2}(\IP^1,\delta_0,\delta_\infty; \rho, \mu)$ the natural projection. 
Therefore 
\[
e_{\bf T}(Ob_\infty)  =(k\sfh)^{-1} (-k\sfh)^{-1}e_{\bf T}(\IE^\vee(k\sfh)) e_{\bf T}(\IE^\vee(-k\sfh))
\]
where $\IE$ is the Hodge bundle. Using Mumford's relation
$c(\IE^\vee) c(\IE) = 1$,
 one further obtains
\be\label{eq:obsB}
\bal 
 e_{\bf T}(Ob_\infty) =  (-1)^{h-1} (k\sfh)^{2h-2}. 
\eal 
\ee
Therefore the rubber factors in \eqref{eq:factorA} reduce to  
\be\label{eq:factorB}
GW^\bullet(h,\rho, \mu)_\sim =  (-1)^{h-1} (k\sfh)^{2h-2}
 \int_{[\om^\bullet_{h_2}
(\IP^1,\delta_0,\delta_\infty; \rho, \mu) ]^{vir}} {1\over -{\sf u}- \psi_0}. 
\ee

Next note that the equation \eqref{eq:factorA} yields a 
gluing formula for the partition function 
\be\label{eq:resfactA}
Z_\mu(Y,\Delta; g_s) = \sum_{h\in \IZ} g_s^{2h-2+l(\mu)} 
GW^\bullet(h,\mu)
\ee
of residual local invariants with relative conditions $\mu$. 
Let 
\[
S_{\mu}(g_s) = \sum_{h\in \IZ} g_s^{2h-2+l(\mu)} GW^\bullet_s(h,\mu)
\]
be the generating function of simple residual invariants with relative conditions and 
\[
R_{\rho, \mu}(g_s) = \sum_{h \in \IZ} g_s^{2h-2+l(\mu)+l(\rho)} 
GW^\bullet(h,\rho, \mu)_\sim
\]
be the generating function of rubber relative invariants with relative conditions $(\rho, \mu)$. 
Then equation \eqref{eq:resfactA} yields 
\be\label{eq:resfactB} 
Z_\mu(Y,\Delta; g_s) = S_{\mu}(g_s) + \sum_{|\rho|=|\lambda|} (ik\sfh)^{2l(\rho)}\zeta(\rho)S_{\rho}(g_s)R_{\rho, \mu}( g_s).
\ee

\subsection{Simple fixed loci}\label{simple}
The domain of a simple fixed map is a union 
\[ 
\Sigma = \Sigma_0 \cup \big(\cup_{j=1}^{l} \Lambda_{i}\big)
\]
where 
\begin{itemize} 
\item $\Sigma_0$ is a possibly disconnected genus $h$ component mapped to the fixed point $p\in C$, and 
\item the components $\Lambda_{1}, \ldots, \Lambda_{l}$ are projective lines attached to $\Sigma_0$, which are 
mapped in a torus invariant fashion to $C$
with some degrees $d_{i}, \ldots, d_{l}\geq 1$. 
\end{itemize} 
The partition $\mu$ encoding the contact conditions at $\delta_\infty$ is determined by $(d_1, \ldots, d_l)$. 

The equivariant residues of simple fixed loci  are evaluated 
by standard localization computations 
in terms of  
Hodge integrals on moduli spaces of curves with marked points. 
Such computations have been done in detail for example in 
\cite{MV_formula,Hodge_unknot, Curves_TQFT,Two_partition},
the results being in agreement with  the open string Gromov-Witten invariants 
computed in \cite{Open_relative,KL}. Omitting the details the end result is a
closed form expression for 
the generating function $S^\bullet_{\alpha}(g_s)$ provided by the 
Marino-Vafa formula
\cite{Framed_knots}, proven in \cite{MV_formula,Hodge_unknot}.  
To write this formula explicitly, 
for any Young diagram $\nu$ let 
$s_\nu(x_1, x_2, \ldots)$ are the corresponding Schur function and let $c(\nu) = \sum_{\Box\in \nu}(a(\Box)-l(\Box))/2$ be the 
content of $\nu$. As above, if $\nu=(1^{k_1}2^{k_2}\cdots)$, let $\zeta(\nu) = \prod_{j\geq 1} k_j! j^{k_j}$. 
Moreover let $\chi^\nu$ denote the character of the irreducible representation of the symmetric group  $\CS_{|\nu|}$ determined by $\nu$. For any Young diagram $\mu$ with 
$|\mu|=|\nu|$  let $\chi^\nu(\mu)$ denote the value of $\chi^\nu$ on the 
conjugacy class determined by $\mu$. 
Then the Marino-Vafa formula reads
\be\label{eq:ressimpleA} 
S_\mu(Y,\dl; g_s) = (ik{\sf u})^{-l(\mu)} \sum_{|\nu|=|\mu|} 
{\chi^{\nu}(\mu) \over \zeta(\mu)} q^{(k+1)c(\nu)}s_\nu(\uq).
\ee
where $q=e^{ig_s}$ and $\uq=(q^{1/2}, q^{3/2}, \ldots)$. 

\subsection{Rubber integrals and Hurwitz numbers}\label{rubberHurwitz}

In order to finish the computation one has to evaluate the rubber integrals in equation \eqref{eq:factorB}, which has been done in \cite{GW_completed, Virasoro_curves, Two_partition}. 
First note that 
\[ 
\bal 
GW^\bullet(h,\rho, \mu)_\sim = (-1)^{h} (k\sfh)^{2h-2}
\sum_{n\geq 0} (-1)^{n} \sfh^{-n-1}\langle \rho,n|\mu\rangle^\sim_h
\eal
\]
where 
\[ 
\langle\rho,n|\mu\rangle^\sim_h = \int_{[\om_{h}^\bullet(\IP^1,\delta_0, \delta_\infty; \rho, \mu)_\sim]^{vir}} \psi_0^n.  
\]
These rubber correlators have been evaluated by rigidification in 
\cite{GW_completed, Virasoro_curves, Two_partition}, the 
 final formulas being expressed in terms of Hurwitz numbers, 
\be\label{eq:rubberinvB}
\langle\rho,n|\mu\rangle^\sim_h = 
{H_h^\bullet({\rho,\mu})\over (n+1)!}
\ee
if $n+1 = 2h-2+l(\rho)+l(\mu)$, respectively $\langle\rho,n|\mu\rangle^\sim_h =0$ for all other values of $n$. 

For completeness, recall that 
the double Hurwitz number $H_h^\bullet(\rho,\mu)$ 
is the  weighted number of genus $h$ disconnected $d:1$ covers of $\IP^1$ with 
fixed branch locus of type \[
\big(\mu, \nu, \underbrace{(2,1^{d-1}), \ldots, 
(2,{1^{d-1}})}_{n(h,\rho,\mu)} \big)
\] 
where $d=|\rho|=|\mu|$, and 
\[
n(h,\rho,\mu) = 2h-2 +l(\rho)+l(\mu).
\]
Each cover is weighted by the inverse of the order of its automorphism 
group. The double Hurwitz numbers are given by the following combinatorial formula \cite{Mirror_elliptic}
\be\label{eq:Hurwitz_I} 
\bal
H^\bullet_h(\rho,\mu) & = \sum_{|\nu|=d} c(\nu)^{n(h,\rho,\mu)}
 {\chi^\nu(\rho)\over \zeta(\rho)}{\chi^\nu(\mu)\over \zeta(\mu)}
\eal
\ee
where $c(\nu)=\sum_{\Box\in \nu} (a(\Box)-l(\Box))$ is the content of $\nu$. 

Using equation \eqref{eq:rubberinvB}, the generating series of rubber invariants is written as 
\[ 
\bal 
R_{\rho,\mu}(g_s) & = (ik\sfh)^{-(l(\rho)+l(\mu))} \sum_{\substack{h\in \IZ\\ 
2h-2+l(\rho)+l(\mu) \geq 1}} (-ikg_s)^{2h-2+l(\rho)+l(\mu)} 
{H_h^\bullet(\rho,\mu)\over (2h-2+l(\rho)+l(\mu))!}\\
& = (ik\sfh)^{-(l(\rho)+l(\mu))}\sum_{n\geq 1} (-ikg_s)^{n} 
{H_{h(n,\rho,\mu)}\over n!} 
\eal 
\]
where 
\[
h(n,\rho,\mu) =\left\{\begin{array}{ll} 
(n-l(\rho)-l(\mu)+2)/2, & {\rm if}\  n-l(\rho)-l(\mu)\
{\rm even} \\
0, & {\rm otherwise.}\end{array}\right.
\]
Next note that formula \eqref{eq:Hurwitz_I} yields the following identity 
\be\label{eq:Hurwitz_II} 
\sum_{n\geq 1} {t^n \over n!}
H^\bullet_{h(n,\rho,\mu)}(\rho,\mu) = 
\sum_{\nu} \big(e^{tc(\nu)}-1\big)
{\chi^\nu(\rho)\over \zeta(\rho)} {\chi^\nu(\mu)\over \zeta(\mu)}.
\ee
for any variable $t$. 
Therefore the final formula for the rubber generating series for two partitions $(\rho, \mu)$ with $|\rho|=|\mu|=d$ is 
\be\label{eq:rubberinvC} 
R_{\rho,\mu}(g_s)  = (ik\sfh)^{-(l(\rho)+l(\mu))} \sum_{|\nu|=d}
\big(e^{-ikc(\nu)}-1\big)
{\chi^\nu(\rho)\over \zeta(\rho)} {\chi^\nu(\mu)\over \zeta(\mu)}.
\ee

The computation of the $(0,-1)$ cap is concluded by substituting 
formulas
\eqref{eq:ressimpleA}, \eqref{eq:rubberinvC} in equation 
\eqref{eq:resfactB}. Using the orthogonality relations 
\[
\sum_{|\rho|=|\mu|=|\nu|} \zeta(\rho)^{-1} \chi^{\mu}(\rho)\chi^{\nu}(\rho) = 
\delta_{\mu,\nu},
\]
it follows that 
\be\label{eq:relcap}
Z_\mu(Y,\Delta; g_s) = (ik\sfh)^{-l(\mu)}\sum_{|\nu|=|\mu|} 
{\chi^\nu(\mu)\over \zeta(\mu)} 
q^{c(\nu)} s_\nu(\uq).
\ee
Using the identity,
\[ 
q^{c(\nu)} s_\nu(\uq) = s_{\nu^t}(\uq)
\]
this formula agrees with the antidiagonal $(0,-1)$ cap computed in \cite{Loc_curves} for torus weight $\sft = k\sfh$. Indeed, applying the TQFT formalism of \cite{Loc_curves}, reviewed in Section \ref{central}, 
the element of ${\bf GW}_d(S^1)$ determined by the
above formula is 
\[ 
\bal 
\sum_{|\mu|=d} (ik\sfh)^{2l(\mu)} \zeta(\mu) 
Z_\mu(Y,\Delta; g_s)e_\mu = \sum_{|\mu|=d}\sum_{|\nu|=d} 
(ik\sfh)^{l(\mu)}
{\chi^\nu(\mu)} 
q^{c(\nu)} s_\nu(\uq)e_\mu.
\eal 
\]
Using the change of basis formula \eqref{eq:invbasis}, 
one finds 
\[ 
\sum_{|\rho|=d}(ik\sfh)^{d}\zeta(\rho)s_{\rho^t}(\uq) v_\rho
\]
in the indempotent basis \eqref{eq:vbasis}. This is an agreement with 
the formula for the coefficients ${\overline \eta}_\rho$ derived on page 38 of \cite{Loc_curves}.

\bibliography{Relative_ref.bib}
 \bibliographystyle{abbrv}

\appendix
\section{Examples}\label{examples}

{\bf Example 1.} $g=1$, $m=2$, $\mu_1=\mu_2=(2,1)$, $n_1=3$, $n_2=4$.  
\[
\bal 
P_{\um,\un}(u,v)=  
{u}^{30}{v}^{30}+2\,{u}^{29}{v}^{30}-2\,{u}^{29}{v}^{29}+4\,{u}^{28}{v}^{30}-6\,{u}^{28}{v}^{29}+5\,{u}^{27}{v}^{30}+2\,{u}^{28}{v}^{28}\\
\mbox{}-14\,{u}^{27}{v}^{29}+7\,{u}^{26}{v}^{30}+9\,{u}^{27}{v}^{28}-22\,{u}^{26}{v}^{29}+8\,{u}^{25}{v}^{30}-2\,{u}^{27}{v}^{27}\\
\mbox{}+23\,{u}^{26}{v}^{28}-32\,{u}^{25}{v}^{29}+10\,{u}^{24}{v}^{30}-10\,{u}^{26}{v}^{27}+44\,{u}^{25}{v}^{28}-40\,{u}^{24}{v}^{29}\\
\mbox{}+11\,{u}^{23}{v}^{30}+2\,{u}^{26}{v}^{26}-28\,{u}^{25}{v}^{27}+68\,{u}^{24}{v}^{28}-50\,{u}^{23}{v}^{29}+11\,{u}^{22}{v}^{30}\\
\mbox{}+10\,{u}^{25}{v}^{26}-60\,{u}^{24}{v}^{27}+94\,{u}^{23}{v}^{28}-58\,{u}^{22}{v}^{29}+10\,{u}^{21}{v}^{30}-2\,{u}^{25}{v}^{25}\\
\mbox{}+31\,{u}^{24}{v}^{26}-100\,{u}^{23}{v}^{27}+119\,{u}^{22}{v}^{28}-58\,{u}^{21}{v}^{29}+9\,{u}^{20}{v}^{30}-10\,{u}^{24}{v}^{25}\\
\mbox{}+69\,{u}^{23}{v}^{26}-148\,{u}^{22}{v}^{27}+144\,{u}^{21}{v}^{28}-54\,{u}^{20}{v}^{29}+6\,{u}^{19}{v}^{30}+2\,{u}^{24}{v}^{24}\\
\mbox{}-32\,{u}^{23}{v}^{25}+124\,{u}^{22}{v}^{26}-194\,{u}^{21}{v}^{27}+144\,{u}^{20}{v}^{28}-48\,{u}^{19}{v}^{29}+4\,{u}^{18}{v}^{30}\\
\mbox{}+10\,{u}^{23}{v}^{24}-74\,{u}^{22}{v}^{25}+190\,{u}^{21}{v}^{26}-240\,{u}^{20}{v}^{27}+136\,{u}^{19}{v}^{28}-34\,{u}^{18}{v}^{29}\\
\mbox{}+{u}^{17}{v}^{30}-2\,{u}^{23}{v}^{23}+32\,{u}^{22}{v}^{24}-140\,{u}^{21}{v}^{25}+261\,{u}^{20}{v}^{26}-240\,{u}^{19}{v}^{27}\\
\mbox{}+118\,{u}^{18}{v}^{28}-20\,{u}^{17}{v}^{29}-10\,{u}^{22}{v}^{23}+77\,{u}^{21}{v}^{24}-222\,{u}^{20}{v}^{25}+324\,{u}^{19}{v}^{26}\\
\mbox{}-226\,{u}^{18}{v}^{27}+86\,{u}^{17}{v}^{28}-4\,{u}^{16}{v}^{29}+2\,{u}^{22}{v}^{22}-32\,{u}^{21}{v}^{23}+149\,{u}^{20}{v}^{24}\\
\mbox{}-316\,{u}^{19}{v}^{25}+323\,{u}^{18}{v}^{26}-190\,{u}^{17}{v}^{27}+44\,{u}^{16}{v}^{28}+10\,{u}^{21}{v}^{22}-78\,{u}^{20}{v}^{23}\\
\mbox{}+246\,{u}^{19}{v}^{24}-392\,{u}^{18}{v}^{25}+297\,{u}^{17}{v}^{26}-134\,{u}^{16}{v}^{27}+6\,{u}^{15}{v}^{28}-2\,{u}^{21}{v}^{21}\\
\mbox{}+32\,{u}^{20}{v}^{22}-154\,{u}^{19}{v}^{23}+357\,{u}^{18}{v}^{24}-384\,{u}^{17}{v}^{25}+241\,{u}^{16}{v}^{26}\\
\mbox{}-56\,{u}^{15}{v}^{27}-10\,{u}^{20}{v}^{21}+78\,{u}^{19}{v}^{22}-262\,{u}^{18}{v}^{23}+443\,{u}^{17}{v}^{24}-340\,{u}^{16}{v}^{25}\\
\mbox{}+152\,{u}^{15}{v}^{26}-4\,{u}^{14}{v}^{27}+2\,{u}^{20}{v}^{20}-32\,{u}^{19}{v}^{21}+157\,{u}^{18}{v}^{22}-386\,{u}^{17}{v}^{23}\\
\mbox{}+418\,{u}^{16}{v}^{24}-260\,{u}^{15}{v}^{25}+44\,{u}^{14}{v}^{26}+10\,{u}^{19}{v}^{20}-78\,{u}^{18}{v}^{21}+271\,{u}^{17}{v}^{22}\\
\mbox{}-474\,{u}^{16}{v}^{23}+354\,{u}^{15}{v}^{24}-134\,{u}^{14}{v}^{25}+{u}^{13}{v}^{26}-2\,{u}^{19}{v}^{19}+32\,{u}^{18}{v}^{20}\\
\mbox{}-158\,{u}^{17}{v}^{21}+405\,{u}^{16}{v}^{22}-428\,{u}^{15}{v}^{23}+241\,{u}^{14}{v}^{24}-20\,{u}^{13}{v}^{25}\\
\mbox{}-10\,{u}^{18}{v}^{19}+78\,{u}^{17}{v}^{20}-276\,{u}^{16}{v}^{21}+484\,{u}^{15}{v}^{22}-340\,{u}^{14}{v}^{23}+86\,{u}^{13}{v}^{24}\\
\mbox{}+2\,{u}^{18}{v}^{18}-32\,{u}^{17}{v}^{19}+158\,{u}^{16}{v}^{20}-412\,{u}^{15}{v}^{21}+418\,{u}^{14}{v}^{22}-190\,{u}^{13}{v}^{23}\\
\mbox{}+4\,{u}^{12}{v}^{24}+10\,{u}^{17}{v}^{18}-78\,{u}^{16}{v}^{19}+278\,{u}^{15}{v}^{20}-474\,{u}^{14}{v}^{21}+297\,{u}^{13}{v}^{22}\\
\mbox{}-34\,{u}^{12}{v}^{23}-2\,{u}^{17}{v}^{17}+32\,{u}^{16}{v}^{18}-158\,{u}^{15}{v}^{19}+405\,{u}^{14}{v}^{20}-384\,{u}^{13}{v}^{21}\\
\mbox{}+118\,{u}^{12}{v}^{22}-10\,{u}^{16}{v}^{17}+78\,{u}^{15}{v}^{18}-276\,{u}^{14}{v}^{19}+443\,{u}^{13}{v}^{20}-226\,{u}^{12}{v}^{21}\\
\mbox{}+6\,{u}^{11}{v}^{22}+2\,{u}^{16}{v}^{16}-32\,{u}^{15}{v}^{17}+158\,{u}^{14}{v}^{18}-386\,{u}^{13}{v}^{19}+323\,{u}^{12}{v}^{20}\\
\mbox{}-48\,{u}^{11}{v}^{21}+10\,{u}^{15}{v}^{16}-78\,{u}^{14}{v}^{17}+271\,{u}^{13}{v}^{18}-392\,{u}^{12}{v}^{19}+136\,{u}^{11}{v}^{20}\\
\eal 
\]
\[
\bal
\mbox{}-2\,{u}^{15}{v}^{15}+32\,{u}^{14}{v}^{16}-158\,{u}^{13}{v}^{17}+357\,{u}^{12}{v}^{18}-240\,{u}^{11}{v}^{19}+9\,{u}^{10}{v}^{20}\\
\mbox{}-10\,{u}^{14}{v}^{15}+78\,{u}^{13}{v}^{16}-262\,{u}^{12}{v}^{17}+324\,{u}^{11}{v}^{18}-54\,{u}^{10}{v}^{19}+2\,{u}^{14}{v}^{14}\\
\mbox{}-32\,{u}^{13}{v}^{15}+157\,{u}^{12}{v}^{16}-316\,{u}^{11}{v}^{17}+144\,{u}^{10}{v}^{18}+10\,{u}^{13}{v}^{14}-78\,{u}^{12}{v}^{15}\\
\mbox{}+246\,{u}^{11}{v}^{16}-240\,{u}^{10}{v}^{17}+10\,{u}^{9}{v}^{18}-2\,{u}^{13}{v}^{13}+32\,{u}^{12}{v}^{14}-154\,{u}^{11}{v}^{15}\\
\mbox{}+261\,{u}^{10}{v}^{16}-58\,{u}^{9}{v}^{17}-10\,{u}^{12}{v}^{13}+78\,{u}^{11}{v}^{14}-222\,{u}^{10}{v}^{15}+144\,{u}^{9}{v}^{16}\\
\mbox{}+2\,{u}^{12}{v}^{12}-32\,{u}^{11}{v}^{13}+149\,{u}^{10}{v}^{14}-194\,{u}^{9}{v}^{15}+11\,{u}^{8}{v}^{16}+10\,{u}^{11}{v}^{12}\\
\mbox{}-78\,{u}^{10}{v}^{13}+190\,{u}^{9}{v}^{14}-58\,{u}^{8}{v}^{15}-2\,{u}^{11}{v}^{11}+32\,{u}^{10}{v}^{12}-140\,{u}^{9}{v}^{13}\\
\mbox{}+119\,{u}^{8}{v}^{14}-10\,{u}^{10}{v}^{11}+77\,{u}^{9}{v}^{12}-148\,{u}^{8}{v}^{13}+11\,{u}^{7}{v}^{14}+2\,{u}^{10}{v}^{10}-32\,{u}^{9}{v}^{11}\\
\mbox{}+124\,{u}^{8}{v}^{12}-50\,{u}^{7}{v}^{13}+10\,{u}^{9}{v}^{10}-74\,{u}^{8}{v}^{11}+94\,{u}^{7}{v}^{12}-2\,{u}^{9}{v}^{9}+32\,{u}^{8}{v}^{10}-100\,{u}^{7}{v}^{11}\\
\mbox{}+10\,{u}^{6}{v}^{12}-10\,{u}^{8}{v}^{9}+69\,{u}^{7}{v}^{10}-40\,{u}^{6}{v}^{11}+2\,{u}^{8}{v}^{8}-32\,{u}^{7}{v}^{9}+68\,{u}^{6}{v}^{10}+10\,{u}^{7}{v}^{8}\\
\mbox{}-60\,{u}^{6}{v}^{9}+8\,{u}^{5}{v}^{10}-2\,{u}^{7}{v}^{7}+31\,{u}^{6}{v}^{8}-32\,{u}^{5}{v}^{9}-10\,{u}^{6}{v}^{7}+44\,{u}^{5}{v}^{8}+2\,{u}^{6}{v}^{6}-28\,{u}^{5}{v}^{7}\\
\mbox{}+7\,{u}^{4}{v}^{8}+10\,{u}^{5}{v}^{6}-22\,{u}^{4}{v}^{7}-2\,{u}^{5}{v}^{5}+23\,{u}^{4}{v}^{6}-10\,{u}^{4}{v}^{5}+5\,{u}^{3}{v}^{6}+2\,{u}^{4}{v}^{4}-14\,{u}^{3}{v}^{5}\\
\mbox{}+9\,{u}^{3}{v}^{4}-2\,{u}^{3}{v}^{3}+4\,{u}^{2}{v}^{4}-6\,{u}^{2}{v}^{3}+2\,{u}^{2}{v}^{2}+2\,u{v}^{2}-2\,uv+1.
\eal
\]

{\bf Example 2.} $g=1$, $m=2$, $\mu_1=(2,1)$, $\mu_2=(1,1,1)$, $n_1=3$, $n_2=4$.  
\[
\bal 
P_{\umu,\un}(u,v) = 
{u}^{38}{v}^{38}+3\,{u}^{37}{v}^{38}-2\,{u}^{37}{v}^{37}+6\,{u}^{36}{v}^{38}-8\,{u}^{36}{v}^{37}+9\,{u}^{35}{v}^{38}+2\,{u}^{36}{v}^{36}\\
\mbox{}-20\,{u}^{35}{v}^{37}+12\,{u}^{34}{v}^{38}+11\,{u}^{35}{v}^{36}-36\,{u}^{34}{v}^{37}+15\,{u}^{33}{v}^{38}-2\,{u}^{35}{v}^{35}\\
\mbox{}+32\,{u}^{34}{v}^{36}-54\,{u}^{33}{v}^{37}+18\,{u}^{32}{v}^{38}-12\,{u}^{34}{v}^{35}+67\,{u}^{33}{v}^{36}-72\,{u}^{32}{v}^{37}\\
\mbox{}+21\,{u}^{31}{v}^{38}+2\,{u}^{34}{v}^{34}-38\,{u}^{33}{v}^{35}+112\,{u}^{32}{v}^{36}-90\,{u}^{31}{v}^{37}+24\,{u}^{30}{v}^{38}\\
\mbox{}+12\,{u}^{33}{v}^{34}-88\,{u}^{32}{v}^{35}+162\,{u}^{31}{v}^{36}-108\,{u}^{30}{v}^{37}+27\,{u}^{29}{v}^{38}-2\,{u}^{33}{v}^{33}\\
\mbox{}+41\,{u}^{32}{v}^{34}-160\,{u}^{31}{v}^{35}+213\,{u}^{30}{v}^{36}-126\,{u}^{29}{v}^{37}+30\,{u}^{28}{v}^{38}-12\,{u}^{32}{v}^{33}\\
\mbox{}+100\,{u}^{31}{v}^{34}-248\,{u}^{30}{v}^{35}+264\,{u}^{29}{v}^{36}-144\,{u}^{28}{v}^{37}+30\,{u}^{27}{v}^{38}\\
\mbox{}+2\,{u}^{32}{v}^{32}-42\,{u}^{31}{v}^{33}+193\,{u}^{30}{v}^{34}-342\,{u}^{29}{v}^{35}+315\,{u}^{28}{v}^{36}-162\,{u}^{27}{v}^{37}\\
\mbox{}+27\,{u}^{26}{v}^{38}+12\,{u}^{31}{v}^{32}-106\,{u}^{30}{v}^{33}+314\,{u}^{29}{v}^{34}-438\,{u}^{28}{v}^{35}+366\,{u}^{27}{v}^{36}\\
\mbox{}-168\,{u}^{26}{v}^{37}+24\,{u}^{25}{v}^{38}-2\,{u}^{31}{v}^{31}+42\,{u}^{30}{v}^{32}-214\,{u}^{29}{v}^{33}+451\,{u}^{28}{v}^{34}\\
\mbox{}-534\,{u}^{27}{v}^{35}+414\,{u}^{26}{v}^{36}-156\,{u}^{25}{v}^{37}+18\,{u}^{24}{v}^{38}-12\,{u}^{30}{v}^{31}+109\,{u}^{29}{v}^{32}\\
\mbox{}-362\,{u}^{28}{v}^{33}+594\,{u}^{27}{v}^{34}-630\,{u}^{26}{v}^{35}+438\,{u}^{25}{v}^{36}-138\,{u}^{24}{v}^{37}\\
\mbox{}+12\,{u}^{23}{v}^{38}+2\,{u}^{30}{v}^{30}-42\,{u}^{29}{v}^{31}+226\,{u}^{28}{v}^{32}-538\,{u}^{27}{v}^{33}+738\,{u}^{26}{v}^{34}\\
\mbox{}-714\,{u}^{25}{v}^{35}+417\,{u}^{24}{v}^{36}-108\,{u}^{23}{v}^{37}+6\,{u}^{22}{v}^{38}+12\,{u}^{29}{v}^{30}-110\,{u}^{28}{v}^{31}\\
\mbox{}+395\,{u}^{27}{v}^{32}-726\,{u}^{26}{v}^{33}+879\,{u}^{25}{v}^{34}-756\,{u}^{24}{v}^{35}+366\,{u}^{23}{v}^{36}\\
\eal
\]
\[
\bal
\mbox{}-72\,{u}^{22}{v}^{37}-2\,{u}^{29}{v}^{29}+42\,{u}^{28}{v}^{30}-232\,{u}^{27}{v}^{31}+604\,{u}^{26}{v}^{32}-918\,{u}^{25}{v}^{33}\\
\mbox{}+996\,{u}^{24}{v}^{34}-726\,{u}^{23}{v}^{35}+291\,{u}^{22}{v}^{36}-36\,{u}^{21}{v}^{37}-12\,{u}^{28}{v}^{29}+110\,{u}^{27}{v}^{30}\\
\mbox{}-416\,{u}^{26}{v}^{31}+835\,{u}^{25}{v}^{32}-1098\,{u}^{24}{v}^{33}+1044\,{u}^{23}{v}^{34}-630\,{u}^{22}{v}^{35}\\
\mbox{}+192\,{u}^{21}{v}^{36}+2\,{u}^{28}{v}^{28}-42\,{u}^{27}{v}^{29}+235\,{u}^{26}{v}^{30}-652\,{u}^{25}{v}^{31}+1073\,{u}^{24}{v}^{32}\\
\mbox{}-1236\,{u}^{23}{v}^{33}+990\,{u}^{22}{v}^{34}-492\,{u}^{21}{v}^{35}+90\,{u}^{20}{v}^{36}+12\,{u}^{27}{v}^{28}\\
\mbox{}-110\,{u}^{26}{v}^{29}+428\,{u}^{25}{v}^{30}-922\,{u}^{24}{v}^{31}+1287\,{u}^{23}{v}^{32}-1278\,{u}^{22}{v}^{33}\\
\mbox{}+843\,{u}^{21}{v}^{34}-312\,{u}^{20}{v}^{35}-2\,{u}^{27}{v}^{27}+42\,{u}^{26}{v}^{28}-236\,{u}^{25}{v}^{29}+685\,{u}^{24}{v}^{30}\\
\mbox{}-1202\,{u}^{23}{v}^{31}+1431\,{u}^{22}{v}^{32}-1182\,{u}^{21}{v}^{33}+627\,{u}^{20}{v}^{34}-120\,{u}^{19}{v}^{35}\\
\mbox{}-12\,{u}^{26}{v}^{27}+110\,{u}^{25}{v}^{28}-434\,{u}^{24}{v}^{29}+988\,{u}^{23}{v}^{30}-1440\,{u}^{22}{v}^{31}\\
\mbox{}+1449\,{u}^{21}{v}^{32}-972\,{u}^{20}{v}^{33}+360\,{u}^{19}{v}^{34}+2\,{u}^{26}{v}^{26}-42\,{u}^{25}{v}^{27}+236\,{u}^{24}{v}^{28}\\
\mbox{}-706\,{u}^{23}{v}^{29}+1305\,{u}^{22}{v}^{30}-1578\,{u}^{21}{v}^{31}+1302\,{u}^{20}{v}^{32}-672\,{u}^{19}{v}^{33}\\
\mbox{}+90\,{u}^{18}{v}^{34}+12\,{u}^{25}{v}^{26}-110\,{u}^{24}{v}^{27}+437\,{u}^{23}{v}^{28}-1036\,{u}^{22}{v}^{29}\\
\mbox{}+1551\,{u}^{21}{v}^{30}-1554\,{u}^{20}{v}^{31}+1014\,{u}^{19}{v}^{32}-312\,{u}^{18}{v}^{33}-2\,{u}^{25}{v}^{25}\\
\mbox{}+42\,{u}^{24}{v}^{26}-236\,{u}^{23}{v}^{27}+718\,{u}^{22}{v}^{28}-1380\,{u}^{21}{v}^{29}+1665\,{u}^{20}{v}^{30}\\
\mbox{}-1344\,{u}^{19}{v}^{31}+627\,{u}^{18}{v}^{32}-36\,{u}^{17}{v}^{33}-12\,{u}^{24}{v}^{25}+110\,{u}^{23}{v}^{26}\\
\mbox{}-438\,{u}^{22}{v}^{27}+1068\,{u}^{21}{v}^{28}-1620\,{u}^{20}{v}^{29}+1590\,{u}^{19}{v}^{30}-972\,{u}^{18}{v}^{31}\\
\mbox{}+192\,{u}^{17}{v}^{32}+2\,{u}^{24}{v}^{24}-42\,{u}^{23}{v}^{25}+236\,{u}^{22}{v}^{26}-724\,{u}^{21}{v}^{27}+1425\,{u}^{20}{v}^{28}\\
\mbox{}-1692\,{u}^{19}{v}^{29}+1302\,{u}^{18}{v}^{30}-492\,{u}^{17}{v}^{31}+6\,{u}^{16}{v}^{32}+12\,{u}^{23}{v}^{24}\\
\mbox{}-110\,{u}^{22}{v}^{25}+438\,{u}^{21}{v}^{26}-1086\,{u}^{20}{v}^{27}+1644\,{u}^{19}{v}^{28}-1554\,{u}^{18}{v}^{29}\\
\mbox{}+843\,{u}^{17}{v}^{30}-72\,{u}^{16}{v}^{31}-2\,{u}^{23}{v}^{23}+42\,{u}^{22}{v}^{24}-236\,{u}^{21}{v}^{25}+727\,{u}^{20}{v}^{26}\\
\mbox{}-1440\,{u}^{19}{v}^{27}+1665\,{u}^{18}{v}^{28}-1182\,{u}^{17}{v}^{29}+291\,{u}^{16}{v}^{30}-12\,{u}^{22}{v}^{23}\\
\mbox{}+110\,{u}^{21}{v}^{24}-438\,{u}^{20}{v}^{25}+1092\,{u}^{19}{v}^{26}-1620\,{u}^{18}{v}^{27}+1449\,{u}^{17}{v}^{28}\\
\mbox{}-630\,{u}^{16}{v}^{29}+12\,{u}^{15}{v}^{30}+2\,{u}^{22}{v}^{22}-42\,{u}^{21}{v}^{23}+236\,{u}^{20}{v}^{24}-728\,{u}^{19}{v}^{25}\\
\mbox{}+1425\,{u}^{18}{v}^{26}-1578\,{u}^{17}{v}^{27}+990\,{u}^{16}{v}^{28}-108\,{u}^{15}{v}^{29}+12\,{u}^{21}{v}^{22}\\
\mbox{}-110\,{u}^{20}{v}^{23}+438\,{u}^{19}{v}^{24}-1086\,{u}^{18}{v}^{25}+1551\,{u}^{17}{v}^{26}-1278\,{u}^{16}{v}^{27}\\
\mbox{}+366\,{u}^{15}{v}^{28}-2\,{u}^{21}{v}^{21}+42\,{u}^{20}{v}^{22}-236\,{u}^{19}{v}^{23}+727\,{u}^{18}{v}^{24}-1380\,{u}^{17}{v}^{25}\\
\mbox{}+1431\,{u}^{16}{v}^{26}-726\,{u}^{15}{v}^{27}+18\,{u}^{14}{v}^{28}-12\,{u}^{20}{v}^{21}+110\,{u}^{19}{v}^{22}\\
\mbox{}-438\,{u}^{18}{v}^{23}+1068\,{u}^{17}{v}^{24}-1440\,{u}^{16}{v}^{25}+1044\,{u}^{15}{v}^{26}-138\,{u}^{14}{v}^{27}\\
\mbox{}+2\,{u}^{20}{v}^{20}-42\,{u}^{19}{v}^{21}+236\,{u}^{18}{v}^{22}-724\,{u}^{17}{v}^{23}+1305\,{u}^{16}{v}^{24}-1236\,{u}^{15}{v}^{25}\\
\eal
\]
\[
\bal
\mbox{}+417\,{u}^{14}{v}^{26}+12\,{u}^{19}{v}^{20}-110\,{u}^{18}{v}^{21}+438\,{u}^{17}{v}^{22}-1036\,{u}^{16}{v}^{23}\\
\mbox{}+1287\,{u}^{15}{v}^{24}-756\,{u}^{14}{v}^{25}+24\,{u}^{13}{v}^{26}-2\,{u}^{19}{v}^{19}+42\,{u}^{18}{v}^{20}-236\,{u}^{17}{v}^{21}\\
\mbox{}+718\,{u}^{16}{v}^{22}-1202\,{u}^{15}{v}^{23}+996\,{u}^{14}{v}^{24}-156\,{u}^{13}{v}^{25}-12\,{u}^{18}{v}^{19}\\
\mbox{}+110\,{u}^{17}{v}^{20}-438\,{u}^{16}{v}^{21}+988\,{u}^{15}{v}^{22}-1098\,{u}^{14}{v}^{23}+438\,{u}^{13}{v}^{24}\\
\mbox{}+2\,{u}^{18}{v}^{18}-42\,{u}^{17}{v}^{19}+236\,{u}^{16}{v}^{20}-706\,{u}^{15}{v}^{21}+1073\,{u}^{14}{v}^{22}-714\,{u}^{13}{v}^{23}\\
\mbox{}+27\,{u}^{12}{v}^{24}+12\,{u}^{17}{v}^{18}-110\,{u}^{16}{v}^{19}+437\,{u}^{15}{v}^{20}-922\,{u}^{14}{v}^{21}+879\,{u}^{13}{v}^{22}\\
\mbox{}-168\,{u}^{12}{v}^{23}-2\,{u}^{17}{v}^{17}+42\,{u}^{16}{v}^{18}-236\,{u}^{15}{v}^{19}+685\,{u}^{14}{v}^{20}-918\,{u}^{13}{v}^{21}\\
\mbox{}+414\,{u}^{12}{v}^{22}-12\,{u}^{16}{v}^{17}+110\,{u}^{15}{v}^{18}-434\,{u}^{14}{v}^{19}+835\,{u}^{13}{v}^{20}\\
\mbox{}-630\,{u}^{12}{v}^{21}+30\,{u}^{11}{v}^{22}+2\,{u}^{16}{v}^{16}-42\,{u}^{15}{v}^{17}+236\,{u}^{14}{v}^{18}-652\,{u}^{13}{v}^{19}\\
\mbox{}+738\,{u}^{12}{v}^{20}-162\,{u}^{11}{v}^{21}+12\,{u}^{15}{v}^{16}-110\,{u}^{14}{v}^{17}+428\,{u}^{13}{v}^{18}\\
\mbox{}-726\,{u}^{12}{v}^{19}+366\,{u}^{11}{v}^{20}-2\,{u}^{15}{v}^{15}+42\,{u}^{14}{v}^{16}-236\,{u}^{13}{v}^{17}+604\,{u}^{12}{v}^{18}\\
\mbox{}-534\,{u}^{11}{v}^{19}+30\,{u}^{10}{v}^{20}-12\,{u}^{14}{v}^{15}+110\,{u}^{13}{v}^{16}-416\,{u}^{12}{v}^{17}+594\,{u}^{11}{v}^{18}\\
\mbox{}-144\,{u}^{10}{v}^{19}+2\,{u}^{14}{v}^{14}-42\,{u}^{13}{v}^{15}+235\,{u}^{12}{v}^{16}-538\,{u}^{11}{v}^{17}+315\,{u}^{10}{v}^{18}\\
\mbox{}+12\,{u}^{13}{v}^{14}-110\,{u}^{12}{v}^{15}+395\,{u}^{11}{v}^{16}-438\,{u}^{10}{v}^{17}+27\,{u}^{9}{v}^{18}-2\,{u}^{13}{v}^{13}\\
\mbox{}+42\,{u}^{12}{v}^{14}-232\,{u}^{11}{v}^{15}+451\,{u}^{10}{v}^{16}-126\,{u}^{9}{v}^{17}-12\,{u}^{12}{v}^{13}+110\,{u}^{11}{v}^{14}\\
\mbox{}-362\,{u}^{10}{v}^{15}+264\,{u}^{9}{v}^{16}+2\,{u}^{12}{v}^{12}-42\,{u}^{11}{v}^{13}+226\,{u}^{10}{v}^{14}-342\,{u}^{9}{v}^{15}\\
\mbox{}+24\,{u}^{8}{v}^{16}+12\,{u}^{11}{v}^{12}-110\,{u}^{10}{v}^{13}+314\,{u}^{9}{v}^{14}-108\,{u}^{8}{v}^{15}-2\,{u}^{11}{v}^{11}\\
\mbox{}+42\,{u}^{10}{v}^{12}-214\,{u}^{9}{v}^{13}+213\,{u}^{8}{v}^{14}-12\,{u}^{10}{v}^{11}+109\,{u}^{9}{v}^{12}-248\,{u}^{8}{v}^{13}\\
\mbox{}+21\,{u}^{7}{v}^{14}+2\,{u}^{10}{v}^{10}-42\,{u}^{9}{v}^{11}+193\,{u}^{8}{v}^{12}-90\,{u}^{7}{v}^{13}+12\,{u}^{9}{v}^{10}-106\,{u}^{8}{v}^{11}\\
\mbox{}+162\,{u}^{7}{v}^{12}-2\,{u}^{9}{v}^{9}+42\,{u}^{8}{v}^{10}-160\,{u}^{7}{v}^{11}+18\,{u}^{6}{v}^{12}-12\,{u}^{8}{v}^{9}+100\,{u}^{7}{v}^{10}\\
\mbox{}-72\,{u}^{6}{v}^{11}+2\,{u}^{8}{v}^{8}-42\,{u}^{7}{v}^{9}+112\,{u}^{6}{v}^{10}+12\,{u}^{7}{v}^{8}-88\,{u}^{6}{v}^{9}+15\,{u}^{5}{v}^{10}-2\,{u}^{7}{v}^{7}\\
\mbox{}+41\,{u}^{6}{v}^{8}-54\,{u}^{5}{v}^{9}-12\,{u}^{6}{v}^{7}+67\,{u}^{5}{v}^{8}+2\,{u}^{6}{v}^{6}-38\,{u}^{5}{v}^{7}+12\,{u}^{4}{v}^{8}+12\,{u}^{5}{v}^{6}-36\,{u}^{4}{v}^{7}\\
\mbox{}-2\,{u}^{5}{v}^{5}+32\,{u}^{4}{v}^{6}-12\,{u}^{4}{v}^{5}+9\,{u}^{3}{v}^{6}+2\,{u}^{4}{v}^{4}-20\,{u}^{3}{v}^{5}+11\,{u}^{3}{v}^{4}-2\,{u}^{3}{v}^{3}+6\,{u}^{2}{v}^{4}\\
\mbox{}-8\,{u}^{2}{v}^{3}+2\,{u}^{2}{v}^{2}+3u {v}^{2}-2u v +1.
\eal 
\]

{\bf Example 3.} $g=1$, $m=2$, $\mu_1=\mu_2=(2,2)$, $n_1=3$, $n_2=4$.  
\[
\bal 
P_{\umu,\un}(u,v)=
{u}^{58}{v}^{58}+2\,{u}^{57}{v}^{58}-2\,{u}^{57}{v}^{57}+6\,{u}^{56}{v}^{58}-6\,{u}^{56}{v}^{57}+9\,{u}^{55}{v}^{58}+2\,{u}^{56}{v}^{56}\\
\mbox{}-18\,{u}^{55}{v}^{57}+17\,{u}^{54}{v}^{58}+9\,{u}^{55}{v}^{56}-36\,{u}^{54}{v}^{57}+22\,{u}^{53}{v}^{58}-2\,{u}^{55}{v}^{55}\\
\mbox{}+28\,{u}^{54}{v}^{56}-68\,{u}^{53}{v}^{57}+34\,{u}^{52}{v}^{58}-10\,{u}^{54}{v}^{55}+68\,{u}^{53}{v}^{56}-108\,{u}^{52}{v}^{57}\\
\mbox{}+41\,{u}^{51}{v}^{58}+2\,{u}^{54}{v}^{54}-34\,{u}^{53}{v}^{55}+136\,{u}^{52}{v}^{56}-164\,{u}^{51}{v}^{57}+57\,{u}^{50}{v}^{58}\\
\mbox{}+10\,{u}^{53}{v}^{54}-90\,{u}^{52}{v}^{55}+245\,{u}^{51}{v}^{56}-228\,{u}^{50}{v}^{57}+66\,{u}^{49}{v}^{58}-2\,{u}^{53}{v}^{53}\\
\mbox{}+37\,{u}^{52}{v}^{54}-196\,{u}^{51}{v}^{55}+391\,{u}^{50}{v}^{56}-308\,{u}^{49}{v}^{57}+86\,{u}^{48}{v}^{58}-10\,{u}^{52}{v}^{53}\\
\mbox{}+102\,{u}^{51}{v}^{54}-378\,{u}^{50}{v}^{55}+591\,{u}^{49}{v}^{56}-396\,{u}^{48}{v}^{57}+97\,{u}^{47}{v}^{58}\\
\mbox{}+2\,{u}^{52}{v}^{52}-38\,{u}^{51}{v}^{53}+238\,{u}^{50}{v}^{54}-648\,{u}^{49}{v}^{55}+826\,{u}^{48}{v}^{56}-500\,{u}^{47}{v}^{57}\\
\mbox{}+121\,{u}^{46}{v}^{58}+10\,{u}^{51}{v}^{52}-108\,{u}^{50}{v}^{53}+480\,{u}^{49}{v}^{54}-1030\,{u}^{48}{v}^{55}\\
\mbox{}+1122\,{u}^{47}{v}^{56}-612\,{u}^{46}{v}^{57}+134\,{u}^{45}{v}^{58}-2\,{u}^{51}{v}^{51}+38\,{u}^{50}{v}^{52}-262\,{u}^{49}{v}^{53}\\
\mbox{}+873\,{u}^{48}{v}^{54}-1520\,{u}^{47}{v}^{55}+1447\,{u}^{46}{v}^{56}-740\,{u}^{45}{v}^{57}+162\,{u}^{44}{v}^{58}\\
\mbox{}-10\,{u}^{50}{v}^{51}+111\,{u}^{49}{v}^{52}-550\,{u}^{48}{v}^{53}+1445\,{u}^{47}{v}^{54}-2136\,{u}^{46}{v}^{55}\\
\mbox{}+1839\,{u}^{45}{v}^{56}-876\,{u}^{44}{v}^{57}+169\,{u}^{43}{v}^{58}+2\,{u}^{50}{v}^{50}-38\,{u}^{49}{v}^{51}+274\,{u}^{48}{v}^{52}\\
\mbox{}-1040\,{u}^{47}{v}^{53}+2235\,{u}^{46}{v}^{54}-2864\,{u}^{45}{v}^{55}+2254\,{u}^{44}{v}^{56}-1028\,{u}^{43}{v}^{57}\\
\mbox{}+183\,{u}^{42}{v}^{58}+10\,{u}^{49}{v}^{50}-112\,{u}^{48}{v}^{51}+594\,{u}^{47}{v}^{52}-1790\,{u}^{46}{v}^{53}\\
\mbox{}+3242\,{u}^{45}{v}^{54}-3720\,{u}^{44}{v}^{55}+2742\,{u}^{43}{v}^{56}-1150\,{u}^{42}{v}^{57}+177\,{u}^{41}{v}^{58}\\
\mbox{}-2\,{u}^{49}{v}^{49}+38\,{u}^{48}{v}^{50}-280\,{u}^{47}{v}^{51}+1152\,{u}^{46}{v}^{52}-2864\,{u}^{45}{v}^{53}\\
\mbox{}+4499\,{u}^{44}{v}^{54}-4688\,{u}^{43}{v}^{55}+3246\,{u}^{42}{v}^{56}-1232\,{u}^{41}{v}^{57}+177\,{u}^{40}{v}^{58}\\
\mbox{}-10\,{u}^{48}{v}^{49}+112\,{u}^{47}{v}^{50}-618\,{u}^{46}{v}^{51}+2049\,{u}^{45}{v}^{52}-4294\,{u}^{44}{v}^{53}\\
\mbox{}+5977\,{u}^{43}{v}^{54}-5784\,{u}^{42}{v}^{55}+3738\,{u}^{41}{v}^{56}-1256\,{u}^{40}{v}^{57}+159\,{u}^{39}{v}^{58}\\
\mbox{}+2\,{u}^{48}{v}^{48}-38\,{u}^{47}{v}^{49}+283\,{u}^{46}{v}^{50}-1224\,{u}^{45}{v}^{51}+3367\,{u}^{44}{v}^{52}\\
\mbox{}-6118\,{u}^{43}{v}^{53}+7717\,{u}^{42}{v}^{54}-6988\,{u}^{41}{v}^{55}+4062\,{u}^{40}{v}^{56}-1232\,{u}^{39}{v}^{57}\\
\mbox{}+144\,{u}^{38}{v}^{58}+10\,{u}^{47}{v}^{48}-112\,{u}^{46}{v}^{49}+630\,{u}^{45}{v}^{50}-2226\,{u}^{44}{v}^{51}\\
\mbox{}+5199\,{u}^{43}{v}^{52}-8336\,{u}^{42}{v}^{53}+9671\,{u}^{41}{v}^{54}-8150\,{u}^{40}{v}^{55}+4243\,{u}^{39}{v}^{56}\\
\mbox{}-1152\,{u}^{38}{v}^{57}+112\,{u}^{37}{v}^{58}-2\,{u}^{47}{v}^{47}+38\,{u}^{46}{v}^{48}-284\,{u}^{45}{v}^{49}+1268\,{u}^{44}{v}^{50}\\
\mbox{}-3746\,{u}^{43}{v}^{51}+7582\,{u}^{42}{v}^{52}-10968\,{u}^{41}{v}^{53}+11885\,{u}^{40}{v}^{54}-9034\,{u}^{39}{v}^{55}\\
\mbox{}+4189\,{u}^{38}{v}^{56}-1020\,{u}^{37}{v}^{57}+88\,{u}^{36}{v}^{58}-10\,{u}^{46}{v}^{47}+112\,{u}^{45}{v}^{48}\\
\mbox{}-636\,{u}^{44}{v}^{49}+2340\,{u}^{43}{v}^{50}-5916\,{u}^{42}{v}^{51}+10577\,{u}^{41}{v}^{52}-14000\,{u}^{40}{v}^{53}\\
\eal 
\]
\[
\bal
\mbox{}+14025\,{u}^{39}{v}^{54}-9528\,{u}^{38}{v}^{55}+3966\,{u}^{37}{v}^{56}-824\,{u}^{36}{v}^{57}+56\,{u}^{35}{v}^{58}\\
\mbox{}+2\,{u}^{46}{v}^{46}-38\,{u}^{45}{v}^{47}+284\,{u}^{44}{v}^{48}-1292\,{u}^{43}{v}^{49}+4015\,{u}^{42}{v}^{50}\\
\mbox{}-8824\,{u}^{41}{v}^{51}+14169\,{u}^{40}{v}^{52}-17426\,{u}^{39}{v}^{53}+15749\,{u}^{38}{v}^{54}-9528\,{u}^{37}{v}^{55}\\
\mbox{}+3502\,{u}^{36}{v}^{56}-624\,{u}^{35}{v}^{57}+34\,{u}^{34}{v}^{58}+10\,{u}^{45}{v}^{46}-112\,{u}^{44}{v}^{47}\\
\mbox{}+639\,{u}^{43}{v}^{48}-2412\,{u}^{42}{v}^{49}+6453\,{u}^{41}{v}^{50}-12554\,{u}^{40}{v}^{51}+18404\,{u}^{39}{v}^{52}\\
\mbox{}-20824\,{u}^{38}{v}^{53}+16718\,{u}^{37}{v}^{54}-9020\,{u}^{36}{v}^{55}+2850\,{u}^{35}{v}^{56}-416\,{u}^{34}{v}^{57}\\
\mbox{}+14\,{u}^{33}{v}^{58}-2\,{u}^{45}{v}^{45}+38\,{u}^{44}{v}^{46}-284\,{u}^{43}{v}^{47}+1304\,{u}^{42}{v}^{48}-4194\,{u}^{41}{v}^{49}\\
\mbox{}+9817\,{u}^{40}{v}^{50}-17134\,{u}^{39}{v}^{51}+23193\,{u}^{38}{v}^{52}-23538\,{u}^{37}{v}^{53}+16823\,{u}^{36}{v}^{54}\\
\mbox{}-7964\,{u}^{35}{v}^{55}+2128\,{u}^{34}{v}^{56}-232\,{u}^{33}{v}^{57}+6\,{u}^{32}{v}^{58}-10\,{u}^{44}{v}^{45}\\
\mbox{}+112\,{u}^{43}{v}^{46}-640\,{u}^{42}{v}^{47}+2456\,{u}^{41}{v}^{48}-6842\,{u}^{40}{v}^{49}+14208\,{u}^{39}{v}^{50}\\
\mbox{}-22596\,{u}^{38}{v}^{51}+28013\,{u}^{37}{v}^{52}-25108\,{u}^{36}{v}^{53}+15875\,{u}^{35}{v}^{54}-6420\,{u}^{34}{v}^{55}\\
\mbox{}+1422\,{u}^{33}{v}^{56}-96\,{u}^{32}{v}^{57}+2\,{u}^{31}{v}^{58}+2\,{u}^{44}{v}^{44}-38\,{u}^{43}{v}^{45}+284\,{u}^{42}{v}^{46}\\
\mbox{}-1310\,{u}^{41}{v}^{47}+4308\,{u}^{40}{v}^{48}-10568\,{u}^{39}{v}^{49}+19733\,{u}^{38}{v}^{50}-28854\,{u}^{37}{v}^{51}\\
\mbox{}+31799\,{u}^{36}{v}^{52}-25232\,{u}^{35}{v}^{53}+13903\,{u}^{34}{v}^{54}-4740\,{u}^{33}{v}^{55}+762\,{u}^{32}{v}^{56}\\
\mbox{}-24\,{u}^{31}{v}^{57}+10\,{u}^{43}{v}^{44}-112\,{u}^{42}{v}^{45}+640\,{u}^{41}{v}^{46}-2480\,{u}^{40}{v}^{47}\\
\mbox{}+7113\,{u}^{39}{v}^{48}-15538\,{u}^{38}{v}^{49}+26397\,{u}^{37}{v}^{50}-35134\,{u}^{36}{v}^{51}+33951\,{u}^{35}{v}^{52}\\
\mbox{}-23668\,{u}^{34}{v}^{53}+11017\,{u}^{33}{v}^{54}-3080\,{u}^{32}{v}^{55}+294\,{u}^{31}{v}^{56}-8\,{u}^{30}{v}^{57}\\
\mbox{}-2\,{u}^{43}{v}^{43}+38\,{u}^{42}{v}^{44}-284\,{u}^{41}{v}^{45}+1313\,{u}^{40}{v}^{46}-4380\,{u}^{39}{v}^{47}\\
\mbox{}+11115\,{u}^{38}{v}^{48}-21900\,{u}^{37}{v}^{49}+34131\,{u}^{36}{v}^{50}-40020\,{u}^{35}{v}^{51}+33933\,{u}^{34}{v}^{52}\\
\mbox{}-20372\,{u}^{33}{v}^{53}+7953\,{u}^{32}{v}^{54}-1592\,{u}^{31}{v}^{55}+42\,{u}^{30}{v}^{56}-10\,{u}^{42}{v}^{43}\\
\mbox{}+112\,{u}^{41}{v}^{44}-640\,{u}^{40}{v}^{45}+2492\,{u}^{39}{v}^{46}-7292\,{u}^{38}{v}^{47}+16565\,{u}^{37}{v}^{48}\\
\mbox{}-29708\,{u}^{36}{v}^{49}+41846\,{u}^{35}{v}^{50}-42580\,{u}^{34}{v}^{51}+31455\,{u}^{33}{v}^{52}-15768\,{u}^{32}{v}^{53}\\
\mbox{}+4931\,{u}^{31}{v}^{54}-536\,{u}^{30}{v}^{55}+12\,{u}^{29}{v}^{56}+2\,{u}^{42}{v}^{42}-38\,{u}^{41}{v}^{43}+284\,{u}^{40}{v}^{44}\\
\mbox{}-1314\,{u}^{39}{v}^{45}+4424\,{u}^{38}{v}^{46}-11506\,{u}^{37}{v}^{47}+23641\,{u}^{36}{v}^{48}-38818\,{u}^{35}{v}^{49}\\
\mbox{}+47736\,{u}^{34}{v}^{50}-42196\,{u}^{33}{v}^{51}+26439\,{u}^{32}{v}^{52}-10952\,{u}^{31}{v}^{53}+2372\,{u}^{30}{v}^{54}\\
\mbox{}-48\,{u}^{29}{v}^{55}+10\,{u}^{41}{v}^{42}-112\,{u}^{40}{v}^{43}+640\,{u}^{39}{v}^{44}-2498\,{u}^{38}{v}^{45}\\
\mbox{}+7406\,{u}^{37}{v}^{46}-17326\,{u}^{36}{v}^{47}+32473\,{u}^{35}{v}^{48}-47886\,{u}^{34}{v}^{49}+50467\,{u}^{33}{v}^{50}\\
\mbox{}-38400\,{u}^{32}{v}^{51}+19799\,{u}^{31}{v}^{52}-6372\,{u}^{30}{v}^{53}+648\,{u}^{29}{v}^{54}-8\,{u}^{28}{v}^{55}\\
\mbox{}-2\,{u}^{41}{v}^{41}+38\,{u}^{40}{v}^{42}-284\,{u}^{39}{v}^{43}+1314\,{u}^{38}{v}^{44}-4448\,{u}^{37}{v}^{45}\\
\eal
\]
\[
\bal
\mbox{}+11777\,{u}^{36}{v}^{46}-25002\,{u}^{35}{v}^{47}+42804\,{u}^{34}{v}^{48}-54550\,{u}^{33}{v}^{49}+49370\,{u}^{32}{v}^{50}\\
\mbox{}-31348\,{u}^{31}{v}^{51}+13023\,{u}^{30}{v}^{52}-2688\,{u}^{29}{v}^{53}+42\,{u}^{28}{v}^{54}-10\,{u}^{40}{v}^{41}\\
\mbox{}+112\,{u}^{39}{v}^{42}-640\,{u}^{38}{v}^{43}+2501\,{u}^{37}{v}^{44}-7478\,{u}^{36}{v}^{45}+17875\,{u}^{35}{v}^{46}\\
\mbox{}-34684\,{u}^{34}{v}^{47}+53025\,{u}^{33}{v}^{48}-57160\,{u}^{32}{v}^{49}+43870\,{u}^{31}{v}^{50}-22448\,{u}^{30}{v}^{51}\\
\mbox{}+6918\,{u}^{29}{v}^{52}-536\,{u}^{28}{v}^{53}+2\,{u}^{27}{v}^{54}+2\,{u}^{40}{v}^{40}-38\,{u}^{39}{v}^{41}+284\,{u}^{38}{v}^{42}\\
\mbox{}-1314\,{u}^{37}{v}^{43}+4460\,{u}^{36}{v}^{44}-11956\,{u}^{35}{v}^{45}+26033\,{u}^{34}{v}^{46}-46046\,{u}^{33}{v}^{47}\\
\mbox{}+60161\,{u}^{32}{v}^{48}-54904\,{u}^{31}{v}^{49}+34540\,{u}^{30}{v}^{50}-13752\,{u}^{29}{v}^{51}+2372\,{u}^{28}{v}^{52}\\
\mbox{}-24\,{u}^{27}{v}^{53}+10\,{u}^{39}{v}^{40}-112\,{u}^{38}{v}^{41}+640\,{u}^{37}{v}^{42}-2502\,{u}^{36}{v}^{43}\\
\mbox{}+7522\,{u}^{35}{v}^{44}-18266\,{u}^{34}{v}^{45}+36388\,{u}^{33}{v}^{46}-57102\,{u}^{32}{v}^{47}+62238\,{u}^{31}{v}^{48}\\
\mbox{}-47392\,{u}^{30}{v}^{49}+23364\,{u}^{29}{v}^{50}-6372\,{u}^{28}{v}^{51}+294\,{u}^{27}{v}^{52}-2\,{u}^{39}{v}^{39}\\
\mbox{}+38\,{u}^{38}{v}^{40}-284\,{u}^{37}{v}^{41}+1314\,{u}^{36}{v}^{42}-4466\,{u}^{35}{v}^{43}+12070\,{u}^{34}{v}^{44}\\
\mbox{}-26784\,{u}^{33}{v}^{45}+48532\,{u}^{32}{v}^{46}-64342\,{u}^{31}{v}^{47}+58388\,{u}^{30}{v}^{48}-35648\,{u}^{29}{v}^{49}\\
\mbox{}+13023\,{u}^{28}{v}^{50}-1592\,{u}^{27}{v}^{51}+6\,{u}^{26}{v}^{52}-10\,{u}^{38}{v}^{39}+112\,{u}^{37}{v}^{40}\\
\mbox{}-640\,{u}^{36}{v}^{41}+2502\,{u}^{35}{v}^{42}-7546\,{u}^{34}{v}^{43}+18536\,{u}^{33}{v}^{44}-37648\,{u}^{32}{v}^{45}\\
\mbox{}+60049\,{u}^{31}{v}^{46}-65386\,{u}^{30}{v}^{47}+48614\,{u}^{29}{v}^{48}-22448\,{u}^{28}{v}^{49}+4931\,{u}^{27}{v}^{50}\\
\mbox{}-96\,{u}^{26}{v}^{51}+2\,{u}^{38}{v}^{38}-38\,{u}^{37}{v}^{39}+284\,{u}^{36}{v}^{40}-1314\,{u}^{35}{v}^{41}+4469\,{u}^{34}{v}^{42}\\
\mbox{}-12142\,{u}^{33}{v}^{43}+27309\,{u}^{32}{v}^{44}-50272\,{u}^{31}{v}^{45}+66922\,{u}^{30}{v}^{46}-59576\,{u}^{29}{v}^{47}\\
\mbox{}+34540\,{u}^{28}{v}^{48}-10952\,{u}^{27}{v}^{49}+762\,{u}^{26}{v}^{50}+10\,{u}^{37}{v}^{38}-112\,{u}^{36}{v}^{39}\\
\mbox{}+640\,{u}^{35}{v}^{40}-2502\,{u}^{34}{v}^{41}+7558\,{u}^{33}{v}^{42}-18712\,{u}^{32}{v}^{43}+38508\,{u}^{31}{v}^{44}\\
\mbox{}-61844\,{u}^{30}{v}^{45}+66448\,{u}^{29}{v}^{46}-47392\,{u}^{28}{v}^{47}+19799\,{u}^{27}{v}^{48}-3080\,{u}^{26}{v}^{49}\\
\mbox{}+14\,{u}^{25}{v}^{50}-2\,{u}^{37}{v}^{37}+38\,{u}^{36}{v}^{38}-284\,{u}^{35}{v}^{39}+1314\,{u}^{34}{v}^{40}-4470\,{u}^{33}{v}^{41}\\
\mbox{}+12186\,{u}^{32}{v}^{42}-27656\,{u}^{31}{v}^{43}+51295\,{u}^{30}{v}^{44}-67792\,{u}^{29}{v}^{45}+58388\,{u}^{28}{v}^{46}\\
\mbox{}-31348\,{u}^{27}{v}^{47}+7953\,{u}^{26}{v}^{48}-232\,{u}^{25}{v}^{49}-10\,{u}^{36}{v}^{37}+112\,{u}^{35}{v}^{38}\\
\mbox{}-640\,{u}^{34}{v}^{39}+2502\,{u}^{33}{v}^{40}-7564\,{u}^{32}{v}^{41}+18820\,{u}^{31}{v}^{42}-39000\,{u}^{30}{v}^{43}\\
\mbox{}+62450\,{u}^{29}{v}^{44}-65386\,{u}^{28}{v}^{45}+43870\,{u}^{27}{v}^{46}-15768\,{u}^{26}{v}^{47}+1422\,{u}^{25}{v}^{48}\\
\mbox{}+2\,{u}^{36}{v}^{36}-38\,{u}^{35}{v}^{37}+284\,{u}^{34}{v}^{38}-1314\,{u}^{33}{v}^{39}+4470\,{u}^{32}{v}^{40}\\
\mbox{}-12210\,{u}^{31}{v}^{41}+27855\,{u}^{30}{v}^{42}-51632\,{u}^{29}{v}^{43}+66922\,{u}^{28}{v}^{44}-54904\,{u}^{27}{v}^{45}\\
\mbox{}+26439\,{u}^{26}{v}^{46}-4740\,{u}^{25}{v}^{47}+34\,{u}^{24}{v}^{48}+10\,{u}^{35}{v}^{36}-112\,{u}^{34}{v}^{37}\\
\mbox{}+640\,{u}^{33}{v}^{38}-2502\,{u}^{32}{v}^{39}+7567\,{u}^{31}{v}^{40}-18880\,{u}^{30}{v}^{41}+39158\,{u}^{29}{v}^{42}\\
\eal
\]
\[
\bal
\mbox{}-61844\,{u}^{28}{v}^{43}+62238\,{u}^{27}{v}^{44}-38400\,{u}^{26}{v}^{45}+11017\,{u}^{25}{v}^{46}-416\,{u}^{24}{v}^{47}\\
\mbox{}-2\,{u}^{35}{v}^{35}+38\,{u}^{34}{v}^{36}-284\,{u}^{33}{v}^{37}+1314\,{u}^{32}{v}^{38}-4470\,{u}^{31}{v}^{39}\\
\mbox{}+12221\,{u}^{30}{v}^{40}-27920\,{u}^{29}{v}^{41}+51295\,{u}^{28}{v}^{42}-64342\,{u}^{27}{v}^{43}+49370\,{u}^{26}{v}^{44}\\
\mbox{}-20372\,{u}^{25}{v}^{45}+2128\,{u}^{24}{v}^{46}-10\,{u}^{34}{v}^{35}+112\,{u}^{33}{v}^{36}-640\,{u}^{32}{v}^{37}\\
\mbox{}+2502\,{u}^{31}{v}^{38}-7568\,{u}^{30}{v}^{39}+18900\,{u}^{29}{v}^{40}-39000\,{u}^{28}{v}^{41}+60049\,{u}^{27}{v}^{42}\\
\mbox{}-57160\,{u}^{26}{v}^{43}+31455\,{u}^{25}{v}^{44}-6420\,{u}^{24}{v}^{45}+56\,{u}^{23}{v}^{46}+2\,{u}^{34}{v}^{34}\\
\mbox{}-38\,{u}^{33}{v}^{35}+284\,{u}^{32}{v}^{36}-1314\,{u}^{31}{v}^{37}+4470\,{u}^{30}{v}^{38}-12224\,{u}^{29}{v}^{39}\\
\mbox{}+27855\,{u}^{28}{v}^{40}-50272\,{u}^{27}{v}^{41}+60161\,{u}^{26}{v}^{42}-42196\,{u}^{25}{v}^{43}+13903\,{u}^{24}{v}^{44}\\
\mbox{}-624\,{u}^{23}{v}^{45}+10\,{u}^{33}{v}^{34}-112\,{u}^{32}{v}^{35}+640\,{u}^{31}{v}^{36}-2502\,{u}^{30}{v}^{37}\\
\mbox{}+7568\,{u}^{29}{v}^{38}-18880\,{u}^{28}{v}^{39}+38508\,{u}^{27}{v}^{40}-57102\,{u}^{26}{v}^{41}+50467\,{u}^{25}{v}^{42}\\
\mbox{}-23668\,{u}^{24}{v}^{43}+2850\,{u}^{23}{v}^{44}-2\,{u}^{33}{v}^{33}+38\,{u}^{32}{v}^{34}-284\,{u}^{31}{v}^{35}\\
\mbox{}+1314\,{u}^{30}{v}^{36}-4470\,{u}^{29}{v}^{37}+12221\,{u}^{28}{v}^{38}-27656\,{u}^{27}{v}^{39}+48532\,{u}^{26}{v}^{40}\\
\mbox{}-54550\,{u}^{25}{v}^{41}+33933\,{u}^{24}{v}^{42}-7964\,{u}^{23}{v}^{43}+88\,{u}^{22}{v}^{44}-10\,{u}^{32}{v}^{33}\\
\mbox{}+112\,{u}^{31}{v}^{34}-640\,{u}^{30}{v}^{35}+2502\,{u}^{29}{v}^{36}-7568\,{u}^{28}{v}^{37}+18820\,{u}^{27}{v}^{38}\\
\mbox{}-37648\,{u}^{26}{v}^{39}+53025\,{u}^{25}{v}^{40}-42580\,{u}^{24}{v}^{41}+15875\,{u}^{23}{v}^{42}-824\,{u}^{22}{v}^{43}\\
\mbox{}+2\,{u}^{32}{v}^{32}-38\,{u}^{31}{v}^{33}+284\,{u}^{30}{v}^{34}-1314\,{u}^{29}{v}^{35}+4470\,{u}^{28}{v}^{36}\\
\mbox{}-12210\,{u}^{27}{v}^{37}+27309\,{u}^{26}{v}^{38}-46046\,{u}^{25}{v}^{39}+47736\,{u}^{24}{v}^{40}-25232\,{u}^{23}{v}^{41}\\
\mbox{}+3502\,{u}^{22}{v}^{42}+10\,{u}^{31}{v}^{32}-112\,{u}^{30}{v}^{33}+640\,{u}^{29}{v}^{34}-2502\,{u}^{28}{v}^{35}\\
\mbox{}+7567\,{u}^{27}{v}^{36}-18712\,{u}^{26}{v}^{37}+36388\,{u}^{25}{v}^{38}-47886\,{u}^{24}{v}^{39}+33951\,{u}^{23}{v}^{40}\\
\mbox{}-9020\,{u}^{22}{v}^{41}+112\,{u}^{21}{v}^{42}-2\,{u}^{31}{v}^{31}+38\,{u}^{30}{v}^{32}-284\,{u}^{29}{v}^{33}+1314\,{u}^{28}{v}^{34}\\
\mbox{}-4470\,{u}^{27}{v}^{35}+12186\,{u}^{26}{v}^{36}-26784\,{u}^{25}{v}^{37}+42804\,{u}^{24}{v}^{38}-40020\,{u}^{23}{v}^{39}\\
\mbox{}+16823\,{u}^{22}{v}^{40}-1020\,{u}^{21}{v}^{41}-10\,{u}^{30}{v}^{31}+112\,{u}^{29}{v}^{32}-640\,{u}^{28}{v}^{33}\\
\mbox{}+2502\,{u}^{27}{v}^{34}-7564\,{u}^{26}{v}^{35}+18536\,{u}^{25}{v}^{36}-34684\,{u}^{24}{v}^{37}+41846\,{u}^{23}{v}^{38}\\
\mbox{}-25108\,{u}^{22}{v}^{39}+3966\,{u}^{21}{v}^{40}+2\,{u}^{30}{v}^{30}-38\,{u}^{29}{v}^{31}+284\,{u}^{28}{v}^{32}\\
\mbox{}-1314\,{u}^{27}{v}^{33}+4470\,{u}^{26}{v}^{34}-12142\,{u}^{25}{v}^{35}+26033\,{u}^{24}{v}^{36}-38818\,{u}^{23}{v}^{37}\\
\mbox{}+31799\,{u}^{22}{v}^{38}-9528\,{u}^{21}{v}^{39}+144\,{u}^{20}{v}^{40}+10\,{u}^{29}{v}^{30}-112\,{u}^{28}{v}^{31}\\
\mbox{}+640\,{u}^{27}{v}^{32}-2502\,{u}^{26}{v}^{33}+7558\,{u}^{25}{v}^{34}-18266\,{u}^{24}{v}^{35}+32473\,{u}^{23}{v}^{36}\\
\mbox{}-35134\,{u}^{22}{v}^{37}+16718\,{u}^{21}{v}^{38}-1152\,{u}^{20}{v}^{39}-2\,{u}^{29}{v}^{29}+38\,{u}^{28}{v}^{30}\\
\mbox{}-284\,{u}^{27}{v}^{31}+1314\,{u}^{26}{v}^{32}-4470\,{u}^{25}{v}^{33}+12070\,{u}^{24}{v}^{34}-25002\,{u}^{23}{v}^{35}\\
\eal 
\]
\[
\bal
\mbox{}+34131\,{u}^{22}{v}^{36}-23538\,{u}^{21}{v}^{37}+4189\,{u}^{20}{v}^{38}-10\,{u}^{28}{v}^{29}+112\,{u}^{27}{v}^{30}\\
\mbox{}-640\,{u}^{26}{v}^{31}+2502\,{u}^{25}{v}^{32}-7546\,{u}^{24}{v}^{33}+17875\,{u}^{23}{v}^{34}-29708\,{u}^{22}{v}^{35}\\
\mbox{}+28013\,{u}^{21}{v}^{36}-9528\,{u}^{20}{v}^{37}+159\,{u}^{19}{v}^{38}+2\,{u}^{28}{v}^{28}-38\,{u}^{27}{v}^{29}\\
\mbox{}+284\,{u}^{26}{v}^{30}-1314\,{u}^{25}{v}^{31}+4469\,{u}^{24}{v}^{32}-11956\,{u}^{23}{v}^{33}+23641\,{u}^{22}{v}^{34}\\
\mbox{}-28854\,{u}^{21}{v}^{35}+15749\,{u}^{20}{v}^{36}-1232\,{u}^{19}{v}^{37}+10\,{u}^{27}{v}^{28}-112\,{u}^{26}{v}^{29}\\
\mbox{}+640\,{u}^{25}{v}^{30}-2502\,{u}^{24}{v}^{31}+7522\,{u}^{23}{v}^{32}-17326\,{u}^{22}{v}^{33}+26397\,{u}^{21}{v}^{34}\\
\mbox{}-20824\,{u}^{20}{v}^{35}+4243\,{u}^{19}{v}^{36}-2\,{u}^{27}{v}^{27}+38\,{u}^{26}{v}^{28}-284\,{u}^{25}{v}^{29}\\
\mbox{}+1314\,{u}^{24}{v}^{30}-4466\,{u}^{23}{v}^{31}+11777\,{u}^{22}{v}^{32}-21900\,{u}^{21}{v}^{33}+23193\,{u}^{20}{v}^{34}\\
\mbox{}-9034\,{u}^{19}{v}^{35}+177\,{u}^{18}{v}^{36}-10\,{u}^{26}{v}^{27}+112\,{u}^{25}{v}^{28}-640\,{u}^{24}{v}^{29}\\
\mbox{}+2502\,{u}^{23}{v}^{30}-7478\,{u}^{22}{v}^{31}+16565\,{u}^{21}{v}^{32}-22596\,{u}^{20}{v}^{33}+14025\,{u}^{19}{v}^{34}\\
\mbox{}-1256\,{u}^{18}{v}^{35}+2\,{u}^{26}{v}^{26}-38\,{u}^{25}{v}^{27}+284\,{u}^{24}{v}^{28}-1314\,{u}^{23}{v}^{29}\\
\mbox{}+4460\,{u}^{22}{v}^{30}-11506\,{u}^{21}{v}^{31}+19733\,{u}^{20}{v}^{32}-17426\,{u}^{19}{v}^{33}+4062\,{u}^{18}{v}^{34}\\
\mbox{}+10\,{u}^{25}{v}^{26}-112\,{u}^{24}{v}^{27}+640\,{u}^{23}{v}^{28}-2502\,{u}^{22}{v}^{29}+7406\,{u}^{21}{v}^{30}\\
\mbox{}-15538\,{u}^{20}{v}^{31}+18404\,{u}^{19}{v}^{32}-8150\,{u}^{18}{v}^{33}+177\,{u}^{17}{v}^{34}-2\,{u}^{25}{v}^{25}\\
\mbox{}+38\,{u}^{24}{v}^{26}-284\,{u}^{23}{v}^{27}+1314\,{u}^{22}{v}^{28}-4448\,{u}^{21}{v}^{29}+11115\,{u}^{20}{v}^{30}\\
\mbox{}-17134\,{u}^{19}{v}^{31}+11885\,{u}^{18}{v}^{32}-1232\,{u}^{17}{v}^{33}-10\,{u}^{24}{v}^{25}+112\,{u}^{23}{v}^{26}\\
\mbox{}-640\,{u}^{22}{v}^{27}+2501\,{u}^{21}{v}^{28}-7292\,{u}^{20}{v}^{29}+14208\,{u}^{19}{v}^{30}-14000\,{u}^{18}{v}^{31}\\
\mbox{}+3738\,{u}^{17}{v}^{32}+2\,{u}^{24}{v}^{24}-38\,{u}^{23}{v}^{25}+284\,{u}^{22}{v}^{26}-1314\,{u}^{21}{v}^{27}\\
\mbox{}+4424\,{u}^{20}{v}^{28}-10568\,{u}^{19}{v}^{29}+14169\,{u}^{18}{v}^{30}-6988\,{u}^{17}{v}^{31}+183\,{u}^{16}{v}^{32}\\
\mbox{}+10\,{u}^{23}{v}^{24}-112\,{u}^{22}{v}^{25}+640\,{u}^{21}{v}^{26}-2498\,{u}^{20}{v}^{27}+7113\,{u}^{19}{v}^{28}\\
\mbox{}-12554\,{u}^{18}{v}^{29}+9671\,{u}^{17}{v}^{30}-1150\,{u}^{16}{v}^{31}-2\,{u}^{23}{v}^{23}+38\,{u}^{22}{v}^{24}\\
\mbox{}-284\,{u}^{21}{v}^{25}+1314\,{u}^{20}{v}^{26}-4380\,{u}^{19}{v}^{27}+9817\,{u}^{18}{v}^{28}-10968\,{u}^{17}{v}^{29}\\
\mbox{}+3246\,{u}^{16}{v}^{30}-10\,{u}^{22}{v}^{23}+112\,{u}^{21}{v}^{24}-640\,{u}^{20}{v}^{25}+2492\,{u}^{19}{v}^{26}\\
\mbox{}-6842\,{u}^{18}{v}^{27}+10577\,{u}^{17}{v}^{28}-5784\,{u}^{16}{v}^{29}+169\,{u}^{15}{v}^{30}+2\,{u}^{22}{v}^{22}\\
\mbox{}-38\,{u}^{21}{v}^{23}+284\,{u}^{20}{v}^{24}-1314\,{u}^{19}{v}^{25}+4308\,{u}^{18}{v}^{26}-8824\,{u}^{17}{v}^{27}\\
\mbox{}+7717\,{u}^{16}{v}^{28}-1028\,{u}^{15}{v}^{29}+10\,{u}^{21}{v}^{22}-112\,{u}^{20}{v}^{23}+640\,{u}^{19}{v}^{24}\\
\mbox{}-2480\,{u}^{18}{v}^{25}+6453\,{u}^{17}{v}^{26}-8336\,{u}^{16}{v}^{27}+2742\,{u}^{15}{v}^{28}-2\,{u}^{21}{v}^{21}\\
\mbox{}+38\,{u}^{20}{v}^{22}-284\,{u}^{19}{v}^{23}+1313\,{u}^{18}{v}^{24}-4194\,{u}^{17}{v}^{25}+7582\,{u}^{16}{v}^{26}\\
\mbox{}-4688\,{u}^{15}{v}^{27}+162\,{u}^{14}{v}^{28}-10\,{u}^{20}{v}^{21}+112\,{u}^{19}{v}^{22}-640\,{u}^{18}{v}^{23}\\
\eal 
\]
\[
\bal
\mbox{}+2456\,{u}^{17}{v}^{24}-5916\,{u}^{16}{v}^{25}+5977\,{u}^{15}{v}^{26}-876\,{u}^{14}{v}^{27}+2\,{u}^{20}{v}^{20}\\
\mbox{}-38\,{u}^{19}{v}^{21}+284\,{u}^{18}{v}^{22}-1310\,{u}^{17}{v}^{23}+4015\,{u}^{16}{v}^{24}-6118\,{u}^{15}{v}^{25}\\
\mbox{}+2254\,{u}^{14}{v}^{26}+10\,{u}^{19}{v}^{20}-112\,{u}^{18}{v}^{21}+640\,{u}^{17}{v}^{22}-2412\,{u}^{16}{v}^{23}\\
\mbox{}+5199\,{u}^{15}{v}^{24}-3720\,{u}^{14}{v}^{25}+134\,{u}^{13}{v}^{26}-2\,{u}^{19}{v}^{19}+38\,{u}^{18}{v}^{20}\\
\mbox{}-284\,{u}^{17}{v}^{21}+1304\,{u}^{16}{v}^{22}-3746\,{u}^{15}{v}^{23}+4499\,{u}^{14}{v}^{24}-740\,{u}^{13}{v}^{25}\\
\mbox{}-10\,{u}^{18}{v}^{19}+112\,{u}^{17}{v}^{20}-640\,{u}^{16}{v}^{21}+2340\,{u}^{15}{v}^{22}-4294\,{u}^{14}{v}^{23}\\
\mbox{}+1839\,{u}^{13}{v}^{24}+2\,{u}^{18}{v}^{18}-38\,{u}^{17}{v}^{19}+284\,{u}^{16}{v}^{20}-1292\,{u}^{15}{v}^{21}\\
\mbox{}+3367\,{u}^{14}{v}^{22}-2864\,{u}^{13}{v}^{23}+121\,{u}^{12}{v}^{24}+10\,{u}^{17}{v}^{18}-112\,{u}^{16}{v}^{19}\\
\mbox{}+639\,{u}^{15}{v}^{20}-2226\,{u}^{14}{v}^{21}+3242\,{u}^{13}{v}^{22}-612\,{u}^{12}{v}^{23}-2\,{u}^{17}{v}^{17}\\
\mbox{}+38\,{u}^{16}{v}^{18}-284\,{u}^{15}{v}^{19}+1268\,{u}^{14}{v}^{20}-2864\,{u}^{13}{v}^{21}+1447\,{u}^{12}{v}^{22}\\
\mbox{}-10\,{u}^{16}{v}^{17}+112\,{u}^{15}{v}^{18}-636\,{u}^{14}{v}^{19}+2049\,{u}^{13}{v}^{20}-2136\,{u}^{12}{v}^{21}\\
\mbox{}+97\,{u}^{11}{v}^{22}+2\,{u}^{16}{v}^{16}-38\,{u}^{15}{v}^{17}+284\,{u}^{14}{v}^{18}-1224\,{u}^{13}{v}^{19}+2235\,{u}^{12}{v}^{20}\\
\mbox{}-500\,{u}^{11}{v}^{21}+10\,{u}^{15}{v}^{16}-112\,{u}^{14}{v}^{17}+630\,{u}^{13}{v}^{18}-1790\,{u}^{12}{v}^{19}\\
\mbox{}+1122\,{u}^{11}{v}^{20}-2\,{u}^{15}{v}^{15}+38\,{u}^{14}{v}^{16}-284\,{u}^{13}{v}^{17}+1152\,{u}^{12}{v}^{18}\\
\mbox{}-1520\,{u}^{11}{v}^{19}+86\,{u}^{10}{v}^{20}-10\,{u}^{14}{v}^{15}+112\,{u}^{13}{v}^{16}-618\,{u}^{12}{v}^{17}\\
\mbox{}+1445\,{u}^{11}{v}^{18}-396\,{u}^{10}{v}^{19}+2\,{u}^{14}{v}^{14}-38\,{u}^{13}{v}^{15}+283\,{u}^{12}{v}^{16}-1040\,{u}^{11}{v}^{17}\\
\mbox{}+826\,{u}^{10}{v}^{18}+10\,{u}^{13}{v}^{14}-112\,{u}^{12}{v}^{15}+594\,{u}^{11}{v}^{16}-1030\,{u}^{10}{v}^{17}\\
\mbox{}+66\,{u}^{9}{v}^{18}-2\,{u}^{13}{v}^{13}+38\,{u}^{12}{v}^{14}-280\,{u}^{11}{v}^{15}+873\,{u}^{10}{v}^{16}-308\,{u}^{9}{v}^{17}\\
\mbox{}-10\,{u}^{12}{v}^{13}+112\,{u}^{11}{v}^{14}-550\,{u}^{10}{v}^{15}+591\,{u}^{9}{v}^{16}+2\,{u}^{12}{v}^{12}-38\,{u}^{11}{v}^{13}\\
\mbox{}+274\,{u}^{10}{v}^{14}-648\,{u}^{9}{v}^{15}+57\,{u}^{8}{v}^{16}+10\,{u}^{11}{v}^{12}-112\,{u}^{10}{v}^{13}+480\,{u}^{9}{v}^{14}\\
\mbox{}-228\,{u}^{8}{v}^{15}-2\,{u}^{11}{v}^{11}+38\,{u}^{10}{v}^{12}-262\,{u}^{9}{v}^{13}+391\,{u}^{8}{v}^{14}-10\,{u}^{10}{v}^{11}\\
\mbox{}+111\,{u}^{9}{v}^{12}-378\,{u}^{8}{v}^{13}+41\,{u}^{7}{v}^{14}+2\,{u}^{10}{v}^{10}-38\,{u}^{9}{v}^{11}+238\,{u}^{8}{v}^{12}-164\,{u}^{7}{v}^{13}\\
\mbox{}+10\,{u}^{9}{v}^{10}-108\,{u}^{8}{v}^{11}+245\,{u}^{7}{v}^{12}-2\,{u}^{9}{v}^{9}+38\,{u}^{8}{v}^{10}-196\,{u}^{7}{v}^{11}+34\,{u}^{6}{v}^{12}\\
\mbox{}-10\,{u}^{8}{v}^{9}+102\,{u}^{7}{v}^{10}-108\,{u}^{6}{v}^{11}+2\,{u}^{8}{v}^{8}-38\,{u}^{7}{v}^{9}+136\,{u}^{6}{v}^{10}+10\,{u}^{7}{v}^{8}-90\,{u}^{6}{v}^{9}\\
\mbox{}+22\,{u}^{5}{v}^{10}-2\,{u}^{7}{v}^{7}+37\,{u}^{6}{v}^{8}-68\,{u}^{5}{v}^{9}-10\,{u}^{6}{v}^{7}+68\,{u}^{5}{v}^{8}+2\,{u}^{6}{v}^{6}-34\,{u}^{5}{v}^{7}+17\,{u}^{4}{v}^{8}\\
\mbox{}+10\,{u}^{5}{v}^{6}-36\,{u}^{4}{v}^{7}-2\,{u}^{5}{v}^{5}+28\,{u}^{4}{v}^{6}-10\,{u}^{4}{v}^{5}+9\,{u}^{3}{v}^{6}+2\,{u}^{4}{v}^{4}-18\,{u}^{3}{v}^{5}+9\,{u}^{3}{v}^{4}\\
\mbox{}-2\,{u}^{3}{v}^{3}+6\,{u}^{2}{v}^{4}-6\,{u}^{2}{v}^{3}+2\,{u}^{2}{v}^{2}+2\,u{v}^{2}-2\,uv+1
\eal 
\]

\end{document}